\documentclass[10pt,preprint]{article}

\usepackage{amsmath, latexsym, amsfonts, amssymb, amsthm}
\usepackage{graphicx, color,hyperref,dsfont,epsfig,caption,wrapfig,subfig}
\usepackage{pstricks}
\usepackage{movie15}
\hypersetup{
    linktoc=page,
    linkcolor=red,          
    citecolor=blue,        
    filecolor=blue,      
    urlcolor=cyan,
    colorlinks=true           
}

\numberwithin{equation}{section}

\newtheorem{theorem}{Theorem}[section]
\newtheorem{lemma}[theorem]{Lemma}
\newtheorem{proposition}[theorem]{Proposition}

\newtheorem{remark}[theorem]{Remark}
\newtheorem{definition}[theorem]{Definition}

\newcounter{conj}
\newtheorem{conjecture}[conj]{Conjecture}

\theoremstyle{definition}

\DeclareMathSymbol{\leqslant}{\mathalpha}{AMSa}{"36} 
\DeclareMathSymbol{\geqslant}{\mathalpha}{AMSa}{"3E} 
\DeclareMathSymbol{\eset}{\mathalpha}{AMSb}{"3F}     
\renewcommand{\leq}{\;\leqslant\;}                   
\renewcommand{\geq}{\;\geqslant\;}                   

\newcommand{\C}{\mathbb{C}}
\newcommand{\D}{\mathbb{D}}
\newcommand{\R}{\mathbb{R}}
\newcommand{\Z}{\mathbb{Z}}
\renewcommand{\H}{\mathbb{H}}

\newcommand{\E}{\mathds{E}}
\renewcommand{\P}{\mathds{P}}
\newcommand{\ind}{\mathds{1}}

\usepackage{xparse}

\DeclareDocumentCommand \Pmp { m m o} {
\IfNoValueTF{#3}
{P_{#1}^{#2}}
{P_{#1}^{#2}\left(#3\right)}
}

\DeclareDocumentCommand \Emp { m m o} {
\IfNoValueTF{#3}
{E_{#1}^{#2}}
{E_{#1}^{#2}\left[#3\right]}
}

\DeclareDocumentCommand \Pbr { m m m m o } {
\IfNoValueTF{#5}
{P_{#1}^{#2\stackrel{#4}{\rightarrow} #3}}
{P_{#1}^{#2\stackrel{#4}{\rightarrow} #3}\left(#5\right)}
}

\DeclareDocumentCommand \Ebr { m m m m o } {
\IfNoValueTF{#5}
{E_{#1}^{#2\stackrel{#4}{\rightarrow} #3}}
{E_{#1}^{#2\stackrel{#4}{\rightarrow} #3}\left[#5\right]}
}

\def\S{\mathbb{S}}
\def\T{\mathbb{T}}

\def\bi{\begin{itemize}}
\def\ei{\end{itemize}}

\def\bnum{\begin{enumerate}}
\def\enum{\end{enumerate}}
\def\<#1{\langle #1 \rangle}

\title{Liouville Quantum Gravity on the complex tori}

\author{ Fran\c{c}ois David \footnote{Institut de Physique Th\'eorique,
CNRS, URA 2306, CEA, IPhT, Gif-sur-Yvette, France.}, R\'emi Rhodes \footnote{Universit{\'e} Paris-Est Marne la Vall\'ee, LAMA, Champs sur Marne, France.},
 Vincent Vargas \footnote{ENS Paris, DMA, 45 rue d'Ulm,  75005 Paris, France.} }

\date{\vspace{-5ex}}

\begin{document}

\maketitle

\begin{abstract}
In this paper, we construct Liouville Quantum Field Theory (LQFT) on the toroidal topology in the spirit of the 1981 seminal work by Polyakov. Our approach follows the construction carried out by the authors together with A. Kupiainen in the case of the Riemann sphere. The difference is here that the moduli space for complex tori is non trivial. Modular  properties of LQFT are thus investigated. This allows us to sum up the LQFT on complex tori over the moduli space,  to   compute the law of the random Liouville modulus, therefore recovering (and extending) formulae obtained by physicists, and make conjectures about the relationship with random planar maps of genus one, eventually weighted by a conformal field theory and conformally embedded onto the torus.
\end{abstract}

\footnotesize

\noindent{\bf Key words or phrases:}  Liouville Quantum Gravity, complex tori, modular covariance, Gaussian multiplicative chaos, KPZ formula, Polyakov formula.


\normalsize


\tableofcontents

\section{Introduction}
Given a two dimensional connected compact Riemann manifold $(M,g)$ without boundary, one introduces the convex  Liouville functional 
defined on maps $X:M\to\R$
\begin{equation}\label{introaction}
S(X,g):= \frac{1}{4\pi}\int_{M}\big(|\partial^{g}X |^2+QR_{g} X  + 4\pi \mu e^{\gamma X  }\big)\,\lambda_{g} 
\end{equation}
where $\partial^{g}$, $R_{g} $ and   $\lambda_{g}$  respectively stand for the gradient, Ricci scalar curvature and volume form  in the metric $g$ and $Q,\mu,\gamma>0$ are constants to be discussed later. If  $Q=\frac{2}{\gamma}$, finding the minimizer $u$ (if exists) of this functional  allows one to uniformize $(M,g)$. Indeed, one defines a new metric by setting $g'=e^{\gamma u}g$ and the new metric $g'$ has constant Ricci scalar curvature $R_{g'}=-2\pi\mu\gamma^2$. These are the foundations of the classical Liouville theory, also known under the name of uniformization of $2d$ Riemann surfaces. 

\medskip
For probabilists, Liouville quantum field theory (LQFT) can be seen as a randomization of the Liouville action. This randomization is made in the spirit of the Feynmann path integral representation of the Brownian motion and applies to surfaces instead of paths for the Brownian motion. Formally, one looks for the construction of a random field $X$ 
with law given heuristically  in terms of a  functional integral 
\begin{equation}\label{pathintegral}
\E[F(X)]=Z^{-1}\int F(X)e^{-S(X,g)}DX 
\end{equation}
where $Z$ is a normalization constant and $DX$ stands for a formal uniform measure on some space of maps  $X:M\to\R$ (In fact, in the case of LQFT, $X$ turns out to live in the space of distributions rather than in the space of of maps).  This formalism describes the law of the log-conformal factor $X$ of a formal random metric of the form $e^{\gamma X}g$ on $M$. Of course, this description is purely formal and giving a mathematical description of this picture is a longstanding problem since the work of Polyakov \cite{Pol}.  It turns out that for the particular values 
\begin{equation}
\gamma\in ]0,2],\quad Q=\frac{2}{\gamma}+\frac{\gamma}{2},
\end{equation}
this field theory is expected to become a {\it Conformal Field Theory} (CFT for short, see \cite{difrancesco,gaw} for a background on this topic). The rigorous construction of such an object has been carried out in \cite{DKRV} in the case of the Riemann sphere and in \cite{huang} in the case of simply connected Riemann structures with boundary. In the case of the torus, there is an important physics literature on the topic (see \cite{kleb,friedan,gupta} for instance and also the recent work \cite{ambjorn} which performs very interesting numerical simulations in relation with planar maps) but up to now no rigorous construction.

The aim of this paper is precisely to construct rigorously this CFT   when the underlying Riemann manifold $(M,g)$ has the topology of the torus. A specific feature of complex tori is the non triviality of the modular space as the conformal structures on the torus can be parameterized by a fundamental domain for the action of ${\rm PSL}_2(\Z)$ on  the upper half-plane $\H$ (see next section for exact definitions). Hence complex tori are usually seen as a good way to probe a theory with respect to its modular dependency:  the reader  may  consult \cite{wang} (and references therein) for recent progresses concerning the classical Liouville equation on tori. The non triviality of the modular space was also our original motivation to study  genus one LQFT as a way to test the robustness of the approach developed in \cite{DKRV} in view of the more ambitious task of constructing LQFT on all compact Riemann surfaces with arbitrary large genus.  We will show that Liouville quantum gravity on tori  naturally enjoys additional modular invariance properties besides more standard CFT properties.  Another important question is about summing all the possible moduli dependent LQFT on tori over the moduli space. In particular, the modulus of the random surface becomes a random variable, the law of which will be computed in section \ref{sec:LQG}.

Let us further mention that there is a growing literature on the study of random planar maps (RPM) with genus higher than one, see \cite{bett1,bett2,bender, Chapuy, miermont} for instance. In particular, it has been shown in this case that RPMs converge to the so-called "Brownian torus", a toroidal variant of the Brownian map (work in progress by Bettinelli and Miermont). Liouville quantum gravity on tori is conjecturally related to these objects: we will formulate some precise mathematical conjectures in section \ref{sec:conj} in the case of genus one RPMs. As a byproduct, we give explicit conjectures on the law of the random modulus of the scaling limit of random planar maps, eventually weighted by a model of statistical physics at its critical point and conformally embedded onto the torus. 

\subsubsection*{Acknowledgements} The authors wish to thank Timothy Budd, Romain Dujardin and Colin Guillarmou  for fruitful discussions.
\section{Background and notations}\label{back}

There are at least two canonical ways to describe all the possible complex structures on tori. 

First, given an element $\tau$ called the {\it modular parameter} in the upper half-plane $\H=\{z\in\C;{\rm Im}(z)>0\}$, we can consider the associate  two-dimensional lattice $\Gamma_\tau=\Z+\tau\Z$ and the complex   structure  on the torus  $\T_\tau:=\C/( \Z+\tau\Z)$ induced by the canonical projection from $\C$ (equipped with its canonical complex structure) to $\T_\tau$.  Two such complex structures are equivalent if their modulus parameters $\tau,\tau'$ are conjugated by an element of the {\it modular group} $\mathcal{M}$
\begin{equation}\label{defdepsi}
\exists \psi\in\mathcal{M},\quad \psi(\tau)=\tau'.
\end{equation}
The {\it modular group} $\mathcal{M}$ is the group  of all linear fractional transformations of the upper half plane $\H$ of the form
\begin{equation}\label{modular}
 \tau \mapsto \frac{a \tau+b}{c \tau+d},\quad a,b,c,d\in\Z,\quad ad-bc=1.
\end{equation}
The modular group is isomorphic to  the projective special linear group  PSL$_2(\Z)$ and is generated by the maps $ \tau \mapsto \psi_1(\tau):= \tau+1$ and $z\mapsto  \psi_2(\tau):=-1/\tau$. One can check that if \eqref{defdepsi} holds with $\psi$ given by \eqref{modular} then the application 
\begin{equation}\label{applicationconforme}
z \mapsto (c\tau+d)z
\end{equation}
is a bijective conformal map between $\T_{\psi(\tau)}$ and $\T_\tau$.

\medskip 
It is then natural to ask for a fundamental domain for this equivalence class of conformal structures. Let us introduce the  {\it moduli space}   $\mathcal{S}$ defined as the quotient $\H/\mathcal{M}$ equipped with the Riemann structure inherited from $\H$ equipped with the Poincar\'e metric via the canonical projection. We call $\lambda_{\mathcal{S}}=d^2z/\mathrm{Im}(z)^2$ the volume form on  $\mathcal{S}$. A fundamental domain for the action of  $\mathcal{M}$ on $\H$ is classically given by
\begin{equation}
\{z\in\H; |{\rm Re}(z)|<\frac{1}{2},\,|z|>1\}.
\end{equation}
The Riemann surface $\mathcal{S}$ is not compact.

\medskip 
Throughout this paper, we will  use  a different though  equivalent  procedure to describe the complex structures of tori based on the fact that complex structures on a Riemann surface are equivalent to  conformal structures. Our approach will be thus  to  consider a fixed surface with the topology of the torus, say the usual torus $\T= \C/\Gamma_{i}$,  and to modify the metric on $\T$ as a function of $\tau$ in order to define different conformal structures indexed by $\tau$. More precisely,  we  denote by $p_\tau:\T\to \T_\tau $ the mapping defined by 
\begin{equation}
\label{ ptau2x}
\forall x=(x_1,x_2)\in \T,\quad p_\tau(x)=x_1+\tau x_2
\end{equation}
and consider on $\T$ the complex structure inherited from this mapping, which we call $\tau$-complex structure. The corresponding conformal structure can be described as the class of all metrics on $\T$ conformally equivalent to the metric $\hat{g}_\tau$ that we construct below. 
 
It will be convenient to introduce first the volume form, Laplace-Beltrami operator and gradient with respect to the standard  metric on $\T$, respectively denoted by $\lambda $, $\triangle $ and $\partial$. We will further denote by $\partial_i$ and $\partial_{ij}$ the partial derivatives in this metric.  

Then we  consider the metric (with $x=(x_1,x_2)\in\T$)
\begin{equation}
\hat{g}_\tau(x)dx^2=|dx_1+\tau dx_2|^2
\end{equation}
Observe that $\hat{g}_\tau$ is ``flat'', i.e. has curvature $0$. The volume form, Laplace-Beltrami operator and gradient with respect to this metric will be denoted respectively by $\lambda_\tau$, $\triangle_\tau$ and $\partial^\tau$ and are given by
\begin{align}
d\lambda_\tau&= {\rm Im}(\tau)\,d\lambda\\
\triangle_\tau&= {\rm Im}(\tau)^{-2}\Big(|\tau|^2\partial_{11}-2{\rm Re}(\tau)\partial_{12}+\partial_{22}\Big)\\
\partial^\tau&= {\rm Im}(\tau)^{-2}\big(|\tau|^2\partial_1-{\rm Re}(\tau)\partial_2\,,\,-{\rm Re}(\tau)\partial_1+ \partial_2\big)& & 
\end{align}
We will further denote by $m_\tau(f)$ the mean value of the function $f$ in the metric $\hat{g}_\tau$, namely
\begin{equation}
m_\tau(f)=\frac{1}{\lambda_\tau(\T)}\int_\T f\,d\lambda_\tau.
\end{equation}
Recall also that a function $f: \T \mapsto \C$ is holomorphic if and only if it solves the following Cauchy-Riemann equations:
\begin{equation}\label{CR}
-\frac{{\rm Re}(\tau)}{{\rm Im}(\tau)} \frac{\partial f}{ \partial x}+ \frac{1}{{\rm Im}(\tau)} \frac{\partial f}{\partial y}= i \frac{\partial f}{ \partial x}
\end{equation}
In particular, the mapping $p_\tau:\T\to \T_\tau$ is a conformal map when $\T$ is equipped with the $\tau$-complex structure. Notice also that 
\begin{equation}
\int_{\T}|\partial^\tau\varphi|^2_\tau\,d\lambda_\tau=\frac{1}{{\rm Im}(\tau)}\int_{\T}|\tau \partial_1\varphi-\partial_2\varphi|^2\,d\lambda.
\end{equation}
Finally, as the Sobolev spaces $H^1(\T,\hat{g}_\tau)$ do not depend on $\tau$, we will simply denote by $H^1(\T)$ the space $H^1(\T,\hat{g}_\tau)$, whatever the choice of the metric  $\hat{g}_\tau$.

\medskip
Finally, given an element $\psi(z)=\frac{az+b}{cz+d}$ in PSL$_2(\Z)$ (with $ad-bc=1$), we associate   the unimodular transformation 
\begin{equation}
\widetilde{\psi}: x=(x_1,x_2)\in (\T,\hat{g}_{\psi(\tau)})\mapsto \widetilde{\psi}(x)=(dx_1+bx_2,cx_1+ax_2)\in (\T,\hat{g}_{ \tau}).
\end{equation}
The reader may easily check with the Cauchy-Riemann equations \eqref{CR} that this mapping  is a biholomorphism. In fact, the map $p_\tau \circ \widetilde{\psi} \circ p_{\psi(\tau)}^{-1}$ is given by \eqref{applicationconforme}. Furthermore $ \widetilde{\psi\circ \varphi}=\widetilde{  \varphi}\circ \widetilde{\psi}$ for every $\psi,\varphi\in\mathcal{M}$. We further set $ \psi^t (z)=\frac{az+c}{bz+d}$ and correspondingly $ \widetilde{\psi}^t (x)=(dx_1+cx_2,bx_1+ax_2)$ if $\psi(z)=\frac{az+b}{cz+d}\in \mathcal{M}$  (we just switch the role of $b,c$).

\begin{remark}
Throughout the paper, we will make use of the Dedekind function  $\eta$  and the theta function of Jacobi $\nu_1$, the basic properties of which are recalled in the appendix \ref{special}.
\end{remark}

\section{Green functions and Gaussian Free Fields}
Here we collect a few facts about the Green functions and Gaussian Free Fields on tori.
\subsection{Green function}\label{sub:green}

The Green function $G_\tau$ on the torus $\T$ equipped with the metric $\hat{g}_\tau$ is the unique solution to
\begin{equation}
-\triangle_\tau G_\tau(w,\cdot)=2\pi(\delta_{w}-\frac{1}{{\rm Im}(\tau)}),\quad\quad  \int_{\T}G(w,\cdot)\,d\lambda_\tau=0 
\end{equation}
with the properties that $G_\tau(x,w)=G_\tau(w,x)$ and $G_\tau$ is smooth in $(x,w)$ except along the diagonal $x=w$.
Because of the translation invariance of $\triangle$, we have $G_\tau(x,w)=G_\tau(x-w,0)$, hence it is customary to call  $G_\tau(x):=G_\tau(x,0)$ the Green function.

\medskip
We can expand the Green function along the eigenvalues of $-(2\pi)^{-1}\triangle_{\tau}$. For this, let us define for $n,m\in\Z$
$$\forall x=(x_1,x_2)\in\T,\quad f_{n,m}(x)={\rm Im}(\tau)^{-1/2} \,e^{2\pi i nx_1+2\pi im x_2}.$$
It is readily seen that the family $(f_{n,m})_{mn,\in\Z}$ is an orthonormal basis of $L^2(\lambda_\tau)$. Furthermore we have
$$-\triangle_{\tau}f_{n,m}=2\pi \big(\frac{2\pi}{{\rm Im}(\tau)^2}|\tau n-m|^2\big)f_{n,m}.$$
Therefore the family $(f_{n,m})_{n,m}$ is a family of eigenfunctions of the operator $-(2\pi)^{-1}\triangle_{\tau}$ with associate eigenvalues $(\frac{2\pi}{{\rm Im}(\tau)^2}|\tau n-m|^2)_{n,m}$. Hence we have the usual eigenfunction expansion of the Green function
\begin{equation}\label{greeneigen}
 G_\tau(x)= \sum_{(n,m)\not=(0,0)}c_{n,m}(\tau)e^{2\pi i nx_1+2\pi im x_2}
\end{equation}
with
\begin{equation}
c_{n,m}(\tau)=\frac{{\rm Im}(\tau)}{2\pi|n\tau-m|^2}.
\end{equation}

The exact shape of the Green function on tori is well known, see for instance  \cite[Ch. II, Th 5.1]{lang} where the computations rely on  some algebraic and geometric tools from elliptic curves. For the sake of completeness, we give a direct elementary proof in  appendix \ref{proofshape} by adapting the proof in \cite[section 3.3]{borwein}:

\begin{proposition}\label{shapegreen}
We have the following exact expression for the Green function
\begin{equation}
\forall  x=(x_1,x_2)\in \T,\quad G_\tau(x)= \pi{\rm Im}(\tau)x_2^2- \ln\Big|\frac{\vartheta_1(x_1+\tau x_2,\tau)}{\eta(\tau)}\Big|.
\end{equation}
\end{proposition}
From this exact form of the Green function, it is readily seen that
$$G_\tau(z,w)=\ln\frac{1}{|z-w|}+O(|z-w|^2)+c_\tau$$ for some constant $c_\tau$ which only depends on $\tau$. 

\medskip 
Now we  describe   the behavior of the Green function under the  action of PSL$_2(\Z)$.  We claim

\begin{proposition}\label{prop:modulargreen}
The Green function $G_\tau$ possesses the following modular invariance
\begin{equation}\label{modulargreen}
\forall \tau\in\H,\forall \psi \in\mathcal{M},\forall x\in\T,\quad G_{\psi(\tau)}(x)=G_\tau(\widetilde{\psi}(x)).
\end{equation}
\end{proposition}

\begin{proposition}\label{fourier:rel}
Consider any real valued $\varphi \in L^2(\T)$ and expand it as a Fourier series (converging in $L^2(\T)$)
\begin{equation}
\varphi(x)=\sum_{(n,m)\in (\Z^2) }\varphi_{n,m} e^{2\pi i nx_1+2\pi im x_2}.
\end{equation}
Then we have the algebraic relation
\begin{equation}
\forall\psi\in \mathcal{M}, \forall(n,m)\in (\Z^2),\quad \varphi_{\widetilde{\psi^{t}}^{-1}(n,m)}=(\varphi\circ \psi)_{n,m}.
\end{equation}
\end{proposition}

\noindent {\it Proof of Propositions \ref{prop:modulargreen} and \ref{fourier:rel}}. Notice that
\begin{equation}\label{rel:cnm}
c_{n,m}(\psi(\tau))=c_{\widetilde{\psi^{t}}^{-1}(n,m)}(\tau),
\end{equation}
from which it is plain to deduce that the Green functions satisfies the relations
\begin{equation}
G_{\psi(\tau)}(x)=G_\tau(\widetilde{\psi} (x)). 
\end{equation}

Let us remark that this representation is indeed compatible with the structure of PSL$_2(\Z)$. More precisely, a  M\"obius transform $\psi(z)$ of $\T$ uniquely determines the coefficients $a,b,c,d$ up to a global multiplicative sign. Therefore there are two reparametrizations $\widetilde{\psi}(x),-\widetilde{\psi}(x) $ of the torus attached to $\psi$. Furthermore, because the Green function is real valued, we have $G_\tau(-\widetilde{\psi}(x))=\overline{G_\tau}(\widetilde{\psi}(x))=G_\tau(\widetilde{\psi}(x))$, hence our claim. The proof of Proposition \ref{fourier:rel} results from the same argument.\qed
 
 \medskip

\subsection{Gaussian Free Fields}
Here we follow the usual strategy to define the continuum Gaussian Free Field given the eigenfunction expansion of the Green function \eqref{greeneigen}. Though the case of tori is not treated in \cite{dubedat}, the reader can easily adapt the arguments in \cite{dubedat} to establish the claims below.

\medskip
The Gaussian Free Field $X_\tau$ with vanishing mean in the metric $\hat{g}_\tau$ is  defined as the sum
\begin{equation}
X_\tau(x)=\sum_{(n,m)\in (\Z^2)^*}\alpha_{n,m}c_{n,m}(\tau)^{1/2}e^{2\pi i nx_1+2\pi im x_2}
\end{equation}
where the sequence $(\alpha_{n,m})_{(n,m)\in (\Z^2)^*}$ has the law of an i.i.d. sequence of standard Gaussian random variables on some probability space $(\Omega,\mathcal{F},\P)$. The series converges in the dual Sobolev space  $H^{-1}(\T)$ (whose topology does not depend on $\tau$). It is obvious to check that the covariance kernel of this Gaussian random distribution is exactly the function $G_\tau$. Furthermore, this series representation is also convenient to see that, $\P$-almost surely, the mapping $\tau\mapsto X_\tau\in H^{-1}(\T)$ is measurable.  From Proposition \ref{prop:modulargreen}, we deduce

\begin{proposition}{\bf (Modular invariance of the GFF)}\label{prop:modularGFF}
The GFF $X_\tau$ possesses the following modular invariance
\begin{equation}\label{modulargreen}
\forall \tau\in\H,\forall \psi \in\mathcal{M}, \quad X_{\psi(\tau)} =X_\tau\circ \widetilde{\psi}\quad \text{in law}.
\end{equation}
\end{proposition}

Let us now recall a few facts on the partition function of the GFF.  The free field partition function  on the torus $\T_\tau$  equipped with the flat metric $\hat{g}_{\tau}$ is given by
 \begin{equation}\label{partGFF}
 Z^{{\rm FF}}(\tau)=\frac{1}{\sqrt{{\rm Im}(\tau)}|\eta(\tau)|^2}.
\end{equation}
The reader may easily follow the sketch of proof given in  \cite{gaw} based on a regularized determinant to establish rigorously this statement.
Because the maps $\psi_1,\psi_2$ generate PSL$_2(\Z)$, one can check with the help of \eqref{rulesdede} that $Z^{{\rm FF}}(\psi(\tau))=Z^{{\rm FF}}(\tau)$ for any $\psi$ in the modular group $\mathcal{M}$. So $Z^{{\rm FF}}(\tau)$ is a modular invariant function. 

\medskip
Now we want to extend the notion of GFF on tori to conformally equivalent metrics. This concept requires some precision in view of the construction of Liouville quantum gravity on tori. First there are infinitely many different complex structures on $\T$, which are equivalent up to  PSL$_2(\Z)$. So we have to consider a log-conformal factor, call it $\varphi_\tau$, that depends on all the possible  complex structures on $\T$ parameterized by $\tau\in\H$ and that is compatible with such an algebraic structure. Second,  having in mind (in the end) to sum over the complex structures on tori, the mapping $\tau\mapsto \varphi_\tau$ must be smooth enough. This motivates the following definition


\begin{definition}{\bf (Modular  log-conformal factor)}
A real-valued mapping $$\varphi: \tau \in\H\to \varphi_\tau\in H^1(\T) $$ is called a modular  log-conformal factor if it is measurable, satisfies
\begin{equation}
\forall \psi\in \mathcal{M},\quad \varphi_{\psi(\tau)}=\varphi_\tau\circ \widetilde{\psi} 
\end{equation}
and the condition that the mapping
\begin{equation}
\tau\in \H\mapsto \int_\T| \partial^\tau\varphi|^2_\tau\,d\lambda_\tau
\end{equation}
is invariant under the action of the modular group.  
\end{definition}


The partition function of the GFF in the metric $e^{\varphi_\tau}\hat{g}_\tau$ will be denoted by $Z^{{\rm FF}}(e^{\varphi_\tau}\hat{g}_\tau)$. Recall that thanks to the conformal invariance of the GFF (with its central charge $c=1$), $Z^{{\rm FF}}$ is nothing but the function
\begin{equation}\label{GFFpart}
Z^{{\rm FF}}(e^{\varphi_\tau}\hat{g}_\tau)=e^{\frac{1}{96\pi}\int_\T| \partial^\tau\varphi_\tau |^2_\tau\,d\lambda_\tau }Z^{{\rm FF}}( \tau ).
\end{equation}

When $\varphi_\tau=0$, we will simply write $Z^{{\rm FF}}(\tau)$ instead of $Z^{{\rm FF}}( \hat{g}_\tau)$. Notice that it is plain to construct modular log-conformal factors in terms of their Fourier expansion:

\begin{proposition}
Assume that we are given a sequence $(\varphi_{n,m})_{n,m\in\Z}$ of measurable functions on the moduli surface  $\mathcal{S}$ such that  
\begin{equation}
\forall \tau\in \mathcal{S},\quad \sum_{n,m}|\varphi_{n,m}(\tau)|^2<+\infty,\quad \overline{\varphi_{n,m}(\tau)}=\varphi_{-n,-m}(\tau).
\end{equation}
Extend each function to the whole half plane by setting 
\begin{equation}
\forall \psi\in \mathcal{M}, \forall \tau \in \mathcal{S},\quad \varphi_{n,m}(\psi(\tau))=\varphi_{\widetilde{\psi^t}^{-1}(n,m)}(\tau)
\end{equation}
Then the following function $\varphi$ is a  modular   log-conformal factor
\begin{equation}
\varphi_\tau(x)=\sum_{n,m}\varphi_{n,m}(\tau)c_{n,m}(\tau)^{1/2}e^{2\pi i n x_1+2\pi i m x_2 }.
\end{equation}
\end{proposition}
 
\noindent {\it Proof.} 
The modular compatibility is established with the help of Proposition \ref{fourier:rel} and \eqref{rel:cnm}. Finally,
\begin{align*}
\int_\T| \partial^\tau\varphi_\tau|^2_\tau\,d\lambda_\tau=&\frac{1}{{\rm Im}(\tau)}\int_\T|\tau \partial_1\varphi_\tau-\partial_2\varphi_\tau|^2 \,d\lambda \\
=&\frac{1}{{\rm Im}(\tau)}\int_\T|2\pi i\sum_{n,m}(\tau n-m)\varphi_{n,m}(\tau)c_{n,m}(\tau)^{1/2}e^{2\pi i n x_1+2\pi i m x_2 }|^2 \,d\lambda\\
=&\frac{4\pi^2}{{\rm Im}(\tau)}\sum_{n,m}|\tau n-m|^2  |\varphi_{n,m}(\tau)|^2c_{n,m}(\tau) \\
=& 2\pi \sum_{n,m}  |\varphi_{n,m}(\tau)|^2 .
\end{align*}
 This latter quantity is finite for all $\tau$ and is invariant under $\mathcal{M}$ by construction.\qed

\subsection{Circle-average approximations of the GFF}
 As the GFF belongs to $\H^{-1}(\T)$ and is not a function, it will be sometimes convenient to work with a regularized version of this process. In our context, it is convenient to work with circle average regularizations of the GFF. The point here is that the circles are expressed in the metric  $\hat{g}_\tau$  and do not look like circles when seen with Euclidean eyes. This reflects our convention to work in the fixed space $\T$;  this is in fact equivalent to working in $\T_{\tau}$ with Euclidean circle average regularizations.
 
 The $\epsilon$-circle centered at $z_0$ in the metric $\hat{g}_\tau$ is the set $\{z_0+c_\tau(\epsilon e^{i\theta});\theta\in [0,2\pi]\}$ where  
 \begin{equation}
c_\tau(z)={\rm Re}(z)-\frac{{\rm Re}(\tau)}{{\rm Im}(\tau)}{\rm Im}(z)+i\frac{{\rm Im}(z)}{{\rm Im}(\tau)}.
\end{equation} 
 
 Denote by $X_{\tau,\epsilon}$ the  $\epsilon$ circle-average regularization of the GFF $X_\tau$
 \begin{equation}
X_{\tau,\epsilon}(x)=\frac{1}{2\pi}\int_0^{2\pi}X_\tau(x+c_\tau(\epsilon e^{i\theta}))\,d\theta.
\end{equation}
 Of course, the above notation is formal as it does not make sense to evaluate the GFF at $x+c_\tau(\epsilon e^{i\theta})$ and then integrate over $\theta$. Yet it makes sense to evaluate the GFF against the uniform measure over this circle (this is standard ) and this can be understood as a way to give sense to the above integral, hence our notation. The reader may check via the Kolmogorov criterion that we can find a modification of $X_{\tau,\epsilon}$ that is both continuous in $\epsilon$ and $x$. The covariance kernel of this Gaussian process is the regularized Green function
 \begin{equation}\label{reggreen}
\E[X_{\tau,\epsilon}(x)X_{\tau,\epsilon}(0)]=:G_{\tau,\epsilon}(x)=\frac{1}{4\pi^2}\int_0^{2\pi}\int_0^{2\pi}G_{\tau }\big(x+c_\tau(\epsilon e^{i \theta})-c_\tau(\epsilon e^{i \theta'})\big)\,d\theta d\theta'.
\end{equation}

 We claim that:
\begin{proposition}
\label{circlegreentorus}
For each fixed  $\tau\in \H$ we have the convergence 
\begin{equation}
\E[X_{\tau,\epsilon} (x)^2]+\ln \epsilon\to \Theta(\tau)
\end{equation}
as $\epsilon\to 0$, uniformly with respect to $x\in \T$, where
\begin{equation}
\label{Thetatau}
\Theta(\tau)=-\ln 2\pi-2\ln|\eta(\tau)|.
\end{equation}
\end{proposition}

\noindent {\it Proof.}  By translational invariance of the law of the GFF, if suffices to establish the convergence  at one fixed point $x$ to get the uniform convergence over $\T$. Thus we choose $x=0$. Recall the expression of the Green function of Proposition \ref{shapegreen} and write 
\begin{equation}\label{greenF}
G_\tau(x)=\pi {\rm Im}(\tau)x_2^2+\ln|\eta(\tau)|+\ln\frac{1}{|p_\tau(x)|}-\ln |F(p_\tau (x))|
\end{equation}
where we have set $$F_\tau(x)=\big|\frac{\nu_1(p_\tau(x),\tau)}{p_\tau(x)}\big|.$$ As the function $\nu_1(x,\tau)$ is analytic at $x=0$ and $\nu_1(0,\tau)=0$, the function $F$     is continuous at $x=0$ and $F_\tau(0)= \partial_x\nu_1(0,\tau)$. Let us mention that $F_\tau(0)\not =0$. Indeed one has 
\begin{equation}\label{dtheta}
\partial_x \vartheta_1(0,\tau)=2\pi \eta(\tau)^3.
\end{equation}
Let us admit for a while this relation. This shows that $1/F(x)$ is well defined and continuous over a neighborhood of $0$.

Therefore, on a neighborhood of $0$ and plugging the  expression \eqref{greenF}  into the regularized expression \eqref{reggreen} evaluated at $x=0$, we can  write (using the Landau notation)
 \begin{align*} 
\E[X_{\tau,\epsilon}(0)^2] =&  \int_0^{2\pi}\int_0^{2\pi}\ln\frac{1}{|p_\tau(c_\tau(\epsilon e^{i \theta}-\epsilon e^{i \theta'}))|}\,d\theta d\theta'+\ln|\eta(\tau)|-\ln |F_\tau(0)|+o(\epsilon).
\end{align*} 
Observe now that $p_\tau(c_\tau(x))=x$ and that $ \int_0^{2\pi}\int_0^{2\pi}\ln\frac{1}{|e^{i \theta}-e^{i \theta'}|}\,d\theta d\theta'=0$. Indeed, notice first that  
$$ \int_0^{2\pi}\int_0^{2\pi}\ln\frac{1}{|e^{i \theta}-e^{i \theta'}|}\,d\theta d\theta'=2\pi  \int_0^{2\pi}\ln\frac{1}{|1-e^{i \theta}|}\,d\theta =2\pi\lim_{r\uparrow 1} \int_0^{2\pi}\ln\frac{1}{|1-re^{i \theta}|}\,d\theta .$$
Then use the series expansion of $\ln(1-x)$ for $|x|<1$ into the relation
$$\int_0^{2\pi}\ln\frac{1}{|1-re^{i \theta}|}\,d\theta=\frac{1}{2}\int_0^{2\pi}\ln\frac{1}{1-re^{i \theta}}\,d\theta+\frac{1}{2}\int_0^{2\pi}\ln\frac{1}{1-re^{-i \theta}}\,d\theta.$$
The relation follows easily. Hence,  we get 
 \begin{align*} 
\E[X_{\tau,\epsilon}(0)^2] =& \ln\frac{1}{\epsilon} +\ln|\eta(\tau)|-\ln |F_\tau(0)|+o(\epsilon).
\end{align*} 
Our claim follows by taking the limit as $\epsilon\to 0$ and by using \eqref{dtheta}.

Let us finally show \eqref{dtheta}. By defining the analytic function
$$\xi(z)=q^{\frac{1}{12}}\prod_{n=1}^{\infty}(1-q^{2n}e^{2i\pi z}),\quad q=e^{i\pi \tau},$$
it is readily seen that from \eqref{releteta} that
$$\vartheta_1(z,\tau)=-i e^{i\pi z}\xi(1)\xi(z)(1-e^{-2i\pi z})\xi(-z).$$
By differentiating at $z=0$ and noticing that $\xi(1)=\eta(\tau)$, it is obvious to check that
$$\partial_z \vartheta_1(0,\tau)=2\pi \eta(\tau)^3, $$
which completes the proof.
\qed

\subsection{Gaussian multiplicative chaos}
Next we turn to the construction of the  interaction measure $ e^{\gamma X}\,d\lambda_{\hat{g}_\tau}$. Recall that $X$ is distribution valued and therefore it is not clear at first sight how to define such a random measure. This context falls under the scope of the theory of Gaussian multiplicative chaos \cite{cf:Kah}. We explain below how we  use it.
 
 \medskip

Define   the random measure
\begin{equation}\label{Meps}
M_{\gamma,\tau,\epsilon}:=\epsilon^{\frac{\gamma^2}{2}}e^{\gamma ( X_{\tau,\epsilon}-\frac{Q}{2}\ln{\rm Im}(\tau))}  \,d\lambda_{\hat{g}_\tau}.
\end{equation}
 
\begin{proposition}\label{law}
For $\gamma\in ]0,2[$,   the  following limit exists  in probability
\begin{equation}
M_{\gamma,\tau}:=\lim_{\epsilon\to 0}M_{\gamma,\tau,\epsilon}
=e^{\frac{\gamma^2}{2}\Theta_{\tau} -\frac{\gamma Q}{2}\ln{\rm Im}(\tau)} \lim_{\epsilon\to 0}e^{\gamma X_{\tau,\epsilon}-\frac{\gamma^2}{2} \E[ X_{\tau,\epsilon}^2] } \,d\lambda_{\hat{g}_\tau}
\end{equation}
in the sense of weak convergence of measures. This limiting measure is almost surely non trivial, without atom and is   (up to the multiplicative constant $e^{\frac{\gamma^2}{2}\Theta_{\tau} -\frac{\gamma Q}{2}\ln{\rm Im}(\tau)}$) a Gaussian multiplicative chaos of the field $X_{\tau}$ with respect to the measure $\lambda_{\hat{g}_\tau}$.

Furthermore, for any modular transformation $\psi\in\mathcal{M}$, we have
\begin{equation}
(X_{\psi(\tau)}, M_{\gamma,\psi(\tau)})\stackrel{law}{=}(X_{\tau}\circ\widetilde{\psi }  , \widetilde{\psi}^{-1}*M_{\gamma,\tau}),
\end{equation}
where $ \widetilde{\psi}^{-1}*M_{\gamma,\tau}$ is the pushforward measure of $M_{\gamma,\tau}$ by $ \widetilde{\psi}^{-1}$.
\end{proposition} 
 
\begin{remark}
Recall the definition of the push forward of a measure: given two measured spaces $(A,\mathcal{A})$ and $(B,\mathcal{B})$, a measurable map $f: (A,\mathcal{A})\to(B,\mathcal{B})$ and a measure $m$ on  $(A,\mathcal{A})$, the push forward $ f*m$ of $f$ by $m$ is the measure on $(B,\mathcal{B})$ defined by $f*m(C)=m(f^{-1}(C))$ for $C\in \mathcal{B}$.
\end{remark} 
\noindent {\it Proof.} The convergence in probability results  from the main result in \cite{shamov}. Then the fact that the measure is non trivial, diffuse and is a Gaussian multiplicative chaos are standard facts of Gaussian multiplicative chaos theory (the reader may consult \cite{review} for instance).

\medskip
What has to be checked is the modular transformation rule. In what follows, in order to avoid confusion, we will further index  with $\tau$ a circle average regularized process  to indicate the metric with respect to which circles have to be understood: for instance $(X_{\tau}\circ\widetilde{\psi })_{\epsilon,\psi(\tau)}$ means that  $X_{\tau}\circ\widetilde{\psi }$ has been regularized at scale $\epsilon$ in the metric $\hat{g}_{\psi(\tau)}$ . From Proposition \ref{prop:modularGFF}, we have
\begin{equation}
(X_{\psi(\tau)}, M_{\gamma,\psi(\tau),\epsilon})\stackrel{law}{=}(X_{\tau}\circ\widetilde{\psi }  , \epsilon^{\frac{\gamma^2}{2}}e^{\gamma ( (X_{\tau}\circ\widetilde{\psi })_{\epsilon,\psi(\tau)}-\frac{Q}{2}\ln{\rm Im}(\psi(\tau)))}  \,d\lambda_{\hat{g}_{\psi(\tau)}}).
\end{equation}
Let us rewrite the latter measure as
\begin{align}
 \epsilon^{\frac{\gamma^2}{2}}&e^{\gamma ( (X_{\tau}\circ\widetilde{\psi })_{\epsilon,\psi(\tau)}-\frac{Q}{2}\ln{\rm Im}(\psi(\tau)))}  \,d\lambda_{\hat{g}_{\psi(\tau)}}\nonumber\\
 =&\frac{{\rm Im}(\psi(\tau))}{{\rm Im}(\tau)}e^{-\frac{\gamma Q}{2}\ln{\rm Im}(\tau)}e^{\frac{\gamma^2}{2}\big(\E[(X_{\tau}\circ\widetilde{\psi })_{\epsilon,\psi(\tau)}^2]+\ln\epsilon\big)}e^{\gamma  (X_{\tau}\circ\widetilde{\psi })_{\epsilon,\psi(\tau)}-\frac{\gamma^2}{2}\E[\gamma ( (X_{\tau}\circ\widetilde{\psi })_{\epsilon,\psi(\tau)}^2]}\,d\lambda_{\hat{g}_{ \tau}}.\label{reformeas}
\end{align}
We  use  the main result of \cite{shamov} to claim that the random measures 
$$
e^{\gamma X_{\tau,\epsilon}\circ \widetilde{ \psi} -\frac{\gamma^2}{2}\E[X_{\tau,\epsilon}\circ \widetilde{ \psi}^2]}\,d\lambda_{\hat{g}_\tau}
$$
and 
$$
e^{\gamma (X_{\tau}\circ\widetilde{ \psi})_{\epsilon,\psi(\tau)} -\frac{\gamma^2}{2}\E[ (X_{\tau}\circ \widetilde{ \psi})_{\epsilon,\psi(\tau)}^2]}\,d\lambda_{\hat{g}_\tau}
$$
converge in probability to the same random measure on $\T$, which is nothing but the Gaussian multiplicative chaos $e^{\gamma X_{\tau }\circ \widetilde{ \psi} -\frac{\gamma^2}{2}\E[X_{\tau }\circ \widetilde{ \psi}^2]}\,d\lambda_{\hat{g}_\tau}$. Furthermore, since  $\widetilde{ \psi}$ maps $\T$ onto itself and has determinant $1$, we have that 
$$e^{\gamma X_{\tau }\circ \widetilde{ \psi} -\frac{\gamma^2}{2}\E[X_{\tau }\circ \widetilde{ \psi}^2]}\,d\lambda_{\hat{g}_\tau}= \widetilde{\psi}^{-1}*e^{\gamma X_{\tau } -\frac{\gamma^2}{2}\E[X_{\tau } ^2]}\,d\lambda_{\hat{g}_\tau}.$$
By using Proposition \ref{circlegreentorus} and taking the limit as $\epsilon\to 0$ in  \eqref{reformeas}, we deduce  that
$$(X_{\psi(\tau)}, M_{\gamma,\psi(\tau)}) \stackrel{law}{=}(X_{\tau}\circ\widetilde{\psi } ,\frac{{\rm Im}(\psi(\tau))}{{\rm Im}(\tau)}e^{-\frac{\gamma Q}{2}\ln{\rm Im}(\psi(\tau))}e^{\frac{\gamma^2}{2}\Theta_{\psi(\tau)}}\widetilde{\psi}^{-1}* M_{\gamma,\tau}). $$
Finally, by using the relations (derive the second one by using \eqref{rulesdede})
\begin{equation}\label{rel:imeta}
{\rm Im}(\psi(\tau))={\rm Im}(\tau)|\psi'(\tau)|\quad \text{ and }\quad |\eta(\psi(\tau))|=|\psi'(\tau)|^{-\frac{1}{4}}|\eta(\tau)|,
\end{equation}
it is readily seen   that
$$\frac{{\rm Im}(\psi(\tau))}{{\rm Im}(\tau)}e^{-\frac{\gamma Q}{2}\ln{\rm Im}(\psi(\tau))}e^{\frac{\gamma^2}{2}\Theta_{\psi(\tau)}}= e^{-\frac{\gamma Q}{2}\ln{\rm Im}( \tau)}e^{\frac{\gamma^2}{2}\Theta_{ \tau}}.$$
This completes the proof.\qed
\section{Liouville Quantum Gravity on tori}
 
 As we have overviewed all the necessary background, we can now investigate the construction of the Liouville quantum field theory on tori. Basically, in what follows, we  construct the LQFT for each complex structure $\tau\in\H$. These LQFT will possess some modular invariance properties as well as standard features of conformal field theories, which we  address in this section.
 
 Let us now fix   $\tau\in\H$: this parameter characterizes the complex structure on $\T$ in terms of the metric $\hat{g}_\tau$ put on $\T$. We further fix any   modular  conformal factor $\varphi_\tau$ and consider the metric  $g_\tau=e^{\varphi_\tau}\hat{g}_\tau$, which is conformally equivalent to $\hat{g}_\tau$. This fixes the metric on $\T$ with respect to which we will construct LQFT. We also introduce the GFF $X_{g_\tau}$ with vanishing mean with respect to $g_\tau$; we denote by $X_{g_\tau,\epsilon}$ the associated $\epsilon$ circle-average regularization. One has the following identity in law 
 \begin{equation}\label{identitylaw}
 X_{g_\tau}-m_{\tau} (X_{g_\tau}) \overset{law}{=} X_\tau.
 \end{equation}  
 
Now we consider the two main parameters of LQFT, namely
\begin{equation}
\gamma\in]0,2[\quad \text{ and }\quad \mu>0,
\end{equation}
which are respectively called coupling constant and cosmological constant. 

Finally, we choose  distinct marked points $(z_i,\alpha_i)_{i=1,\dots,n}$ with $z_i\in\T$ and $\alpha_i\in\R$ for all $i=1,\dots,n$. 
Recall that one marked point serves to remove any degree of freedom with respect to the automorphisms of the tori (three marked points are needed for the sphere or the disk, see  \cite{DKRV,huang}).  $n>1$ marked points allows to define correlation functions for  Liouville QFT and LQG on the torus.

\subsection{Definition of the partition function}
We give now the definition function of the partition function in the metric $g_\tau$ for some fixed $\tau\in\H$. It is expressed in terms of a limiting procedure so that we first   define the $\epsilon$-regularized partition function for every measurable bounded function $F$ of $H^{-1}(\T)$
 \begin{align}\label{actiontori2}
 \Pi_{\gamma, \mu}^{(z_i,\alpha_i)_i }  ( g_\tau ,F,\epsilon) 
 :=& Z^{{\rm FF}}( g_\tau) \int_\R  \E\Big[F( c+  X_{g_\tau,\epsilon}  -\frac{Q}{2}\ln{\rm Im}(\tau)+\frac{Q}{2} \varphi_\tau   )  \prod_i \epsilon^{\frac{\alpha_i^2}{2}}e^{\alpha_i (c+  X_{g_\tau,\epsilon} (z_i) -\frac{Q}{2}\ln{\rm Im}(\tau)+\frac{Q}{2} \varphi_\tau(z_i))} \\
 &  \exp\Big( -\frac{Q}{4\pi}\int_{\T }R_{g_\tau}(c+X_{g_\tau}  )\,\lambda_{g_\tau} - \mu \epsilon^{\frac{\gamma^2}{2}}\int_{\T} e^{\gamma ( c+X_{g_\tau,\epsilon}-\frac{Q}{2}\ln{\rm Im}(\tau)+\frac{Q}{2} \varphi_\tau)}  \,d\lambda_{\hat{g}_\tau}   \Big) \Big]\,dc,\nonumber
\end{align}
We want to inquire if the limit  
\begin{align}
 \Pi_{\gamma, \mu}^{(z_i,\alpha_i)_i }  ( g_\tau ,F):=& \lim_{\epsilon\to 0} \Pi_{\gamma, \mu}^{(z_i,\alpha_i)_i }  ( g_\tau ,F,\epsilon),
\end{align}
 exists and is non trivial, meaning that it belongs to $]0,+\infty[$. Non-triviality of the limit will be expressed in terms of the following two conditions, which correspond to  the \textbf{Seiberg bounds} for the torus
 \begin{align}\label{seiberg1}
 &\sum_i\alpha_i> 0,\\ 
 &\forall i,\quad \alpha_i<Q\label{seiberg2}.
 \end{align}
 
The first observation  we make is that  existence and non triviality of the limit do  not depend on the modular  log-conformal factor $\varphi_\tau$. Moreover we have an explicit relation between the cases $\varphi_\tau=0$ and $\varphi_\tau\not=0$:

\begin{theorem}{\bf (Weyl anomaly)}\label{th:weyl}\\
Let us consider a bounded measurable functional
$$F: h \in  H^{-1}(\T) \to F(h) \in\R.$$   Then the two following statements are equivalent:\\
1) the limit $ \Pi_{\gamma, \mu}^{(z_i,\alpha_i)_i }  ( g_\tau ,F)$ exists and is non trivial, \\ 
2) the limit $ \Pi_{\gamma, \mu}^{(z_i,\alpha_i)_i }  ( \hat{g}_\tau ,F)$ exists and is non trivial.

In case they both exist and are  non trivial, we have the following Weyl anomaly
\begin{equation}\label{Weylformula}
 \ln\frac{\Pi_{\gamma, \mu}^{(z_i,\alpha_i)_i } \big(g_\tau,F \big)}{ \Pi_{\gamma, \mu}^{(z_i,\alpha_i)_i }  \big(\hat{g}_\tau ,F \big)}  = \frac{1+6Q^2}{96\pi} \int_{\T}|\partial^\tau   \varphi_\tau|_\tau^2 \,d\lambda_{\hat{g}_\tau} .
 \end{equation}
In particular, this shows that the {\bf central charge} of LQFT on tori is $c_L=1+6Q^2$.
\end{theorem}

\begin{remark} This property is an expected feature of all the conformal field theories (CFT) that can be rigorously constructed. Informally speaking, if $Z_{{\rm CFT}}(g)$ stands for the partition function of any CFT in the background metric $g$ on a compact surface $S$ without boundary, then it is expected that for all conformally equivalent metric $g'=e^{\varphi}g$, we have the relation
\begin{equation}
Z_{{\rm CFT}}(e^{\varphi}g)=e^{\frac{c}{96\pi} \int_{S}|\partial   \varphi|_g^2 +2R_g\varphi \,d\lambda_g }Z_{{\rm CFT}}(g).
\end{equation} 
Therefore, for any CFT,  we should have an explicit dependence on the conformal factor of the metric in terms of the Liouville action (with vanishing cosmological constant) $\frac{1}{96\pi}\int_{S}|\partial   \varphi|_g^2 +2R_g\varphi \,d\lambda_g$ up to some multiplicative  factor $c$,  called the central charge.
\end{remark}
\noindent {\it Proof.}  

First, we notice that one can replace $X_{g_\tau}$ by $X_\tau$ in the definition \eqref{actiontori2} by making the change of variable $c'=  c+m_{\tau} (X_{g_\tau})$ and using identity \eqref{identitylaw}. Therefore, in the sequel, we will work with $X_\tau$ in place of $X_{g_\tau}$.  

Now, the strategy is to use the Girsanov transform in the $\epsilon$- regularized expression \eqref{actiontori2} to the exponential term 
$$\exp\Big( - \frac{Q}{4\pi}\int_{\T}R_{g_\tau}  X_{ \tau } \,d\lambda_{g_\tau}    \Big),$$
which has the effect of shifting the field $X_\tau$ by $$-  \frac{Q}{4\pi}\int_{\T}R_{g_\tau}  G_{\tau}(\cdot-z) \,\lambda_{g_\tau}(dz) =-  \frac{Q}{2}(\varphi_\tau-m_{\tau}(\varphi_\tau)).$$
This Girsanov transform is valid provided that we have renormalized the expression \eqref{actiontori2} by the exponential of 
$\frac{1}{2}$ of the variance of the term $\frac{Q}{4\pi}\int_{\T}R_{g_\tau}  X_{ \tau } \,d\lambda_{g_\tau}  $. This variance  can be computed with the help of the relation for curvature
$$R_{g_\tau}=-e^{-\varphi_\tau }\triangle_\tau\varphi_\tau $$
and is given by
$$\frac{Q^2}{16\pi}\int_{\T}|\partial^\tau   \varphi_\tau|_\tau^2 \,d\lambda_{\hat{g}_\tau} .$$
Then we make the changes  of variables $v=c+\frac{Q}{2}m_{\tau}(\varphi_\tau)$ to get
\begin{align*} 
 \Pi_{\gamma, \mu}^{(z_i,\alpha_i)_i }&   \big(g_\tau,F,\epsilon\big)\\
 =&e^{\frac{6 Q^2}{96\pi}\int_{\T}|\partial^\tau   \varphi_\tau|_\tau^2 \,d\lambda_{\hat{g}_\tau}  } Z^{{\rm FF}}( g_\tau)\int_\R   \E\Big[F( v+X_{  \tau } -\frac{Q}{2}\ln{\rm Im}(\tau)  )\\
 &\prod_i \epsilon^{\frac{\alpha_i^2}{2}}e^{\alpha_i  (v+ X_{\tau,\epsilon} (z_i)-\frac{Q}{2}\ln{\rm Im}(\tau))}   \exp\Big(  - \mu e^{\gamma v}\epsilon^{\frac{\gamma^2}{2}}\int_{\T}e^{\gamma ( X_{\tau,\epsilon} -\frac{Q}{2}\ln{\rm Im}(\tau))}\,d\lambda_{\hat{g}_\tau}    \Big) \Big]\,dv.
 \end{align*}  
 Finally, by using  \eqref{GFFpart} we get
 \begin{align*} 
 \Pi_{\gamma, \mu}^{(z_i,\alpha_i)_i }   \big(g_\tau,F ,\epsilon\big)=&e^{\frac{1+6 Q^2}{96\pi}\int_{\T}|\partial^\tau   \varphi_\tau|_\tau^2 \,d\lambda_{\hat{g}_\tau} }  \Pi_{\gamma, \mu}^{(z_i,\alpha_i)_i }   \big( \hat{g}_\tau ,F,\epsilon\big)\nonumber.
 \end{align*}   
 All our claims are then a straightforward consequence of this relation. \qed 
 
 \medskip
 Notice that the above theorem does not give any information on the existence of the limit. It just allows us to argue that the analysis of the case $\varphi_\tau=0$ is equivalent in some sense to the analysis of the case $\varphi_\tau\not =0$. It is then more convenient to focus on the case   $\varphi_\tau=0$ as the curvature term \eqref{actiontori2} vanishes because of the relation $R_{\hat{g}_\tau}=0$. This observation will be used throughout the remaining part of the paper: the properties of LQFT in the metric $g_\tau$ can be obtained from those of LQFT with metric $\hat{g}_\tau$.   The theorem below investigates the convergence of the partition function \eqref{actiontori2} as $\epsilon\to 0$ in the background metric $\hat{g}_\tau$.

\begin{theorem}{\bf (Convergence of the partition function $\Pi_{\gamma, \mu}^{(z_i,\alpha_i)_i }   ( \hat{g}_\tau ,1)$)}\label{th:seibergtau} \\
1) If \eqref{seiberg1} fails to hold then 
$$\forall \epsilon>0,\quad \Pi_{\gamma, \mu}^{(z_i,\alpha_i)_i }   ( \hat{g}_\tau ,1,\epsilon)=+\infty.$$

\medskip
\noindent 2) Otherwise, i.e. if \eqref{seiberg1} holds,   the limit 
$$\lim_{\epsilon\to 0}\Pi_{\gamma, \mu}^{(z_i,\alpha_i)_i }   ( \hat{g}_\tau ,1,\epsilon):=\Pi_{\gamma,\mu}^{(z_i\alpha_i)_i}  ( \hat{g}_\tau ,1)$$
exists. It is nonzero if and only if \eqref{seiberg2} holds. \\
\medskip
\noindent 3) If \eqref{seiberg1}$+$\eqref{seiberg2} hold, $ \Pi_{\gamma, \mu}^{(z_i,\alpha_i)_i }   ( \hat{g}_\tau ,F)$ is well defined for all bounded $F:\R\times\H^{-1}(\T)\to \R$ and is given by
$$\Pi_{\gamma, \mu}^{(z_i,\alpha_i)_i }   ( \hat{g}_\tau ,F)= Z^{{\rm FF}}( \hat{g}_\tau)e^{C_\tau(\mathbf{z})} \int_\R e^{c\sum_i\alpha_i} \E\Big[F(c+ X_{\tau}+H_{\tau} -\frac{Q}{2}\ln{\rm Im}(\tau) )
 \exp\Big(   - \mu e^{\gamma c}\int_{\T}e^{\gamma H_\tau}\,dM_{\gamma,\tau}    \Big) \Big]\,dc$$
where (denoting $\mathbf{z}=\{z_i\}_i$)
$$C_\tau(\mathbf{z})=\sum_{i<j}\alpha_i\alpha_jG_\tau(z_i-z_j)+\frac{\Theta(\tau)}{2}\sum_i\alpha_i^2-\frac{Q}{2}\ln{\rm Im}(\tau)\sum_i\alpha_i \quad \text{and }\quad H^{\mathbf{z}}_\tau(x)=\sum_i\alpha_iG_\tau(x-z_i).$$
\end{theorem} 
 
\noindent {\it Proof. }   use the same strategy as in \cite{DKRV} .\qed
 
\begin{remark}
By combining Theorems \ref{th:weyl} and \ref{th:seibergtau}, we see that \eqref{seiberg1}$+$\eqref{seiberg2} are also necessary and sufficient in order that the partition function in any background metric $g_\tau=e^{\varphi_\tau}\hat{g}_\tau$ (i.e.  the $\epsilon\to 0$ limit of \eqref{actiontori2}) exists and is non trivial.
\end{remark} 

We denote by $M_+(\T)$ the space of positive measures on $\T$ equipped with the topology of weak convergence. More generally when \eqref{seiberg1}+\eqref{seiberg2} hold then we will define $ \Pi_{\gamma, \mu}^{(z_i,\alpha_i)_i }   ( g_\tau ,F)$ for all bounded $F:\H^{-1}(\T)  \times M_+(\T)  \to \R$ by a limit of the form \eqref{actiontori2} where one simply replaces $F( c+  X_{g_\tau,\epsilon} (z_i) -\frac{Q}{2}\ln{\rm Im}(\tau)+\frac{Q}{2} \varphi_\tau(z_i)    )$ in the definition of \eqref{actiontori2} by
\begin{equation*}
F( c+  X_{g_\tau,\epsilon}  -\frac{Q}{2}\ln{\rm Im}(\tau)+\frac{Q}{2} \varphi_\tau   ,  \epsilon^{\frac{\gamma^2}{2}}  e^{\gamma ( c+X_{g_\tau,\epsilon}-\frac{Q}{2}\ln{\rm Im}(\tau)+\frac{Q}{2} \varphi_\tau)}  d\lambda_{\hat{g}_\tau}   ).
\end{equation*}  
Under this more general framework, the Wey anomaly formula \eqref{Weylformula} is still valid and we have the following formula 
\begin{align}
& \Pi_{\gamma, \mu}^{(z_i,\alpha_i)_i }   ( \hat{g}_\tau ,F)  \nonumber \\
& = Z^{{\rm FF}}( \hat{g}_\tau)e^{C_\tau(\mathbf{z})} \int_\R e^{c\sum_i\alpha_i} \E\Big[F(c+ X_{\tau}+H_{\tau}-  \frac{Q}{2}\ln{\rm Im}( \tau ),e^{\gamma c} e^{\gamma H_\tau}\,dM_{\gamma,\tau}  )
 \exp\Big(   - \mu e^{\gamma c}\int_{\T}e^{\gamma H_\tau}\,dM_{\gamma,\tau}    \Big) \Big]\,dc   \label{defpartitiongeneral}  
 \end{align}

\subsection{Main properties}
Here we study the main properties of the LQFT on tori. Once again, we only consider the case  when the  background metric is  $\hat{g}_\tau$. The more general case of metrics that are conformally equivalent to this metric can easily be recovered thanks to Theorem \ref{th:weyl}.

Let us start with the scaling properties with respect to the cosmological constant $\mu$.

\begin{theorem}{\bf (KPZ scaling laws)}\label{KPZscaling}
We have the following exact scaling relation for the Liouville partition function  with insertions $(z_i,\alpha_i)_i$
$$\Pi_{\gamma,\mu}^{(z_i\alpha_i)_i}  (\hat{g}_\tau,1) =\mu^{\frac{-\sum_i\alpha_i}{\gamma}}\Pi_{\gamma,1}^{(z_i\alpha_i)_i}  (\hat{g}_\tau,1).$$ 
\end{theorem}
The proof of this statement is similar to \cite[Th 3.4]{DKRV} and we let the reader check the details (just make the changes of variables $c'=c+\frac{1}{\gamma}\ln \mu$ in the expression for the partition function in Theorem \ref{th:seibergtau}). Notice that the exponent of $\mu$ differs from the case of the sphere and this is mainly due to the fact that tori have vanishing curvature.

\medskip
Now we focus on  modular invariance/covariance of the partition function, namely  we  describe   the behavior of the partition function under the  action of $\mathcal{M}=$PSL$_2(\Z)$). Recall that we have  associated  to each element of the modular group $\psi(z)=\frac{az+b}{cz+d}$ (with $ad-bc=1$)   the conformal map
$$\widetilde{\psi}: x=(x_1,x_2)\in (\T,\hat{g}_{\psi(\tau)} )\mapsto \widetilde{\psi}(x)=(dx_1+bx_2,cx_1+ax_2)\in (\T,\hat{g}_{\tau} ).$$

\begin{theorem}{\bf (Modular invariance and  KPZ formula)}\label{KPZtori}
Consider an element $\psi$ of the modular group $\mathcal{M}$. Then
$$\Pi_{\gamma,\mu}^{(z_i,\alpha_i)_i} (\hat{g}_{\psi(\tau)},1)=\Big(\prod_i|\psi'(\tau)|^{- \triangle_{\alpha_i}}\Big)\Pi_{\gamma,\mu}^{(\widetilde{\psi}(z_i),\alpha_i)_i} (\hat{g}_\tau,1),$$
where the conformal weight is defined by $\triangle_\alpha=\frac{\alpha}{2}(Q-\frac{\alpha}{2})$.
\end{theorem}

Now we introduce the Liouville field $\phi_\tau$ and the Liouville measure $e^{\gamma \phi_\tau}\,d\lambda_{\hat{g}_\tau}$. Basically these two objects respectively  stand  for the log-conformal factor and volume form of some formal random metric $e^{\gamma \phi_\tau}\hat{g}_\tau$ describing the random geometry of interest in LQFT. We will establish  their modular invariance.   

Let $\varphi_\tau$ be a  modular log-conformal factor and set $g_\tau=e^{\varphi_\tau}\hat{g}_\tau$.   We define a probability law $\P^{  \gamma,\mu}_{(z_i,\alpha_i)_i,g_\tau}$ on $ H^{-1}(\T) \times M_+(\T)$ through its functional   expectations
\begin{equation}\label{LoiLiouville}
\E^{  \gamma,\mu}_{(z_i,\alpha_i)_i,g_\tau}[F]=\frac{\Pi_{\gamma,\mu}^{(z_i\alpha_i)_i} (g_\tau,F)}{\Pi_{\gamma,\mu}^{(z_i\alpha_i)_i} (g_\tau,1)},
\end{equation}
defined for all bounded measurable $F$ on $H^{-1}(\T) \times M_+(\T)$.

\begin{definition}{\bf (modular Liouville field/measure)}\label{def:LFtori}\\
The couple of random variables $(\phi_\tau,e^{\gamma \phi_\tau}\,d\lambda_{g_\tau})$, respectively the Liouville field and Liouville measure, taking values in $\H^{-1}(\T)\times M_+(\T)$ is defined under the probability measure $\P^{  \gamma,\mu}_{(z_i,\alpha_i)_i,g_\tau}$ by relation \eqref{LoiLiouville}.

\end{definition}
%
  
From    Theorem \ref{th:weyl} (or rather its extension described at the end of the previous subsection), we see  that  
\begin{proposition}
The law of the couple Liouville field/measure $(\phi_\tau,e^{\gamma \phi_\tau}\,d\lambda_{g_\tau})$ under $\P^{  \gamma,\mu}_{(z_i,\alpha_i)_i,g_\tau}$ does not depend on the metric $g_\tau=e^{\varphi_\tau}\hat{g}_\tau$  in the conformal equivalence class of $\hat{g}_\tau$, i.e. it does not depend on the modular   log-conformal factor $\varphi_\tau$.
\end{proposition}

Now we study the behavior of this law  under modular transformations. Recall that $f*m$ stands for the pushforward of the measure $m$ by the mapping $f$.
\begin{theorem}\label{coroCRtori}
The law of the couple Liouville field/measure is  modular invariant: i.e. the law of $$(\phi_{\psi(\tau)},e^{\gamma \phi_{\psi(\tau)}}\,d\lambda_{\hat{g}_{\psi(\tau)}})\quad \text{ under }\P^{  \gamma,\mu}_{( z_i,\alpha_i)_i,\hat{g}_\tau}$$ is the same as that of 
$$( \phi_{\tau}\circ \widetilde{\psi} +Q\ln|\widetilde{\psi}'(\tau)|,\widetilde{\psi}^{-1}*e^{\gamma \phi_{\tau}}\,d\lambda_{g_{\tau}})\quad \text{ under }\P^{  \gamma,\mu}_{(\widetilde{\psi}(z_i),\alpha_i)_i,\hat{g}_{\tau}}$$ for any modular transformation  $\psi\in\mathcal{M}$.
 \end{theorem} 

\noindent {\it Proof of Theorems \ref{KPZtori} and   \ref{coroCRtori}.} 
Let $F$ be a bounded continuous functional on $\H^{-1}(\T)\times M_+(\T)$. 
From Theorem \ref{th:seibergtau}, we have
\begin{align*}
\E^{  \gamma,\mu}_{( z_i,\alpha_i)_i,\hat{g}_{\psi(\tau)}}&[F(\phi_{\psi(\tau)},e^{\gamma \phi_{\psi(\tau)}}\,d\lambda_{g_{\psi(\tau)}})]\\
=& Z^{{\rm FF}}( \hat{g}_{\psi(\tau)})e^{C_{\psi(\tau)}( \mathbf{z})} \int_\R e^{c\sum_i\alpha_i} \E\Big[F\big(c+X_{\psi(\tau)}+H^{\mathbf{z}}_{\psi(\tau)} -\frac{Q}{2}\ln{\rm Im}( \psi(\tau)) ,e^{\gamma  H^{\mathbf{z}}_{\psi(\tau)} }e^{\gamma c}M_{\gamma,\psi(\tau)}\big)\\
&
 \exp\Big(   - \mu e^{\gamma c}\int_{\T}e^{\gamma H^{ \mathbf{z}}_{\psi(\tau)} }\,dM_{\gamma,\psi(\tau)}    \Big) \Big]\,dc\times \frac{1}{\Pi_{\gamma,\mu}^{( z_i,\alpha_i)_i} (\hat{g}_{\psi(\tau)},1)}.
 \end{align*}
From Proposition \ref{prop:modulargreen}, we have
\begin{align*} 
 C_{\psi(\tau)}(\mathbf{z})=&\sum_{i<j}\alpha_i\alpha_jG_{\psi(\tau)}(z_i-z_j)+\frac{\Theta(\psi(\tau))}{2}\sum_i\alpha_i^2-\frac{Q}{2}\ln{\rm Im}(\psi(\tau))\sum_i\alpha_i\\
 =& \sum_{i<j}\alpha_i\alpha_jG_{ \tau}(\widetilde{\psi}(z_i)-\widetilde{\psi}(z_j))+\frac{\Theta(\psi(\tau))}{2}\sum_i\alpha_i^2-\frac{Q}{2}\ln{\rm Im}(\psi(\tau))\sum_i\alpha_i\\
 =& C_{\tau}(\widetilde{\psi}(\mathbf{z}))+ \ln|\psi'(\tau)|\sum_i(\frac{\alpha_i^2}{4}-\frac{\alpha_i Q}{2}),
\end{align*} 
 where we have used the following relations 
\begin{equation}\label{CImod}
{\rm Im}(\psi(\tau))={\rm Im}(\tau) |\psi'(\tau)|\quad \text{and}\quad \Theta(\psi(\tau))=\Theta(\tau)+\frac{1}{2}\ln|\psi'(\tau)| .
\end{equation}
The first relation results from a direct computation whereas the second is derived   from   \ref{rulesdede} together with the fact that the mappings $\tau\mapsto \tau+1$ and $\tau\mapsto -1/\tau$ generate $\mathcal{M}$.
Furthermore, from Proposition \ref{prop:modulargreen} again
$$ H^{\mathbf{z}}_{\psi(\tau)}(x) =\sum_i\alpha_iG_{\psi(\tau)}(x-z_i)=\sum_i\alpha_iG_{\tau}(\widetilde{\psi}(x)-\widetilde{\psi}(z_i))= H^{\widetilde{\psi}(\mathbf{z})}_{\tau}(\widetilde{\psi}(x)).$$ 
Hence 
\begin{align*} 
 Z^{{\rm FF}}(\hat{g}_{\psi(\tau)})&e^{C_{\psi(\tau)}( \mathbf{z})} \int_\R e^{c\sum_i\alpha_i} \E\Big[F\big(c+X_{\psi(\tau)}+H^{\mathbf{z}}_{\psi(\tau)} -\frac{Q}{2}\ln{\rm Im}( \psi(\tau)) ,e^{\gamma  H^{\mathbf{z}}_{\psi(\tau)} }e^{\gamma c}M_{\gamma,\psi(\tau)}\big)\\
& \exp\Big(   - \mu e^{\gamma c}\int_{\T}e^{\gamma H^{ \mathbf{z}}_{\psi(\tau)} }\,dM_{\gamma,\psi(\tau)}    \Big) \Big]\,dc\\
=&Z^{{\rm FF}}(\hat{g}_{\tau})e^{ C_{\tau}(\widetilde{\psi}(\mathbf{z}))} \int_\R e^{c\sum_i\alpha_i} \E\Big[F\big(c+X_{\psi(\tau)}+H^{\widetilde{\psi}(\mathbf{z})}_{\tau}\circ \widetilde{\psi}
 -\frac{Q}{2}\ln{\rm Im}( \psi(\tau)) ,e^{\gamma  H^{\widetilde{\psi}(\mathbf{z})}_{\tau}\circ \widetilde{\psi}
 }e^{\gamma c}M_{\gamma,\psi(\tau)}\big)\\
& \exp\Big(   - \mu e^{\gamma c}\int_{\T}e^{\gamma H^{\widetilde{\psi}(\mathbf{z})}_{\tau}\circ\widetilde{\psi} }\,dM_{\gamma,\psi(\tau)}    \Big) \Big]\,dc\times \prod_i|\psi'(\tau)|^{-\triangle_{\alpha_i}}.
 \end{align*}
 It remains to use Proposition \ref{law} and \eqref{CImod} again to see that this latter quantity is equal to  
 \begin{align*} 
 \int_\R & e^{c\sum_i\alpha_i} \E\Big[F\big(c+X_{\tau}\circ \widetilde{\psi}+H^{\widetilde{\psi}(\mathbf{z})}_{\tau}(\widetilde{\psi}(x)) -\frac{Q}{2}\ln{\rm Im}(\tau)+Q\ln|\psi'(\tau)| ,\widetilde{\psi}^{-1}*e^{\gamma  H^{\widetilde{\psi}(\mathbf{z})}_{\tau}}e^{\gamma c}M_{\gamma, \tau}\big)\\
& \exp\Big(   - \mu e^{\gamma c}\int_{\T}e^{\gamma H^{ \widetilde{\psi}(\mathbf{z})}_{\tau} }\,dM_{\gamma,\tau}    \Big) \Big]\,dc\times Z^{{\rm FF}} (\hat{g}_{\tau})e^{ C_{\tau}(\widetilde{\psi}(\mathbf{z}))}\prod_i|\psi'(\tau)|^{-\triangle_{\alpha_i}}\\
=&\E^{  \gamma,\mu}_{( \widetilde{\psi}(z_i),\alpha_i)_i,\hat{g}_{\tau}} [F(\phi_\tau\circ \widetilde{\psi}+Q\ln|\psi'(\tau)|,\widetilde{\psi}^{-1}*e^{\gamma \phi_\tau}\,d\lambda_{g_\tau})]\times  \Pi_{\gamma,\mu}^{( \widetilde{\psi}(z_i),\alpha_i)_i} (\hat{g}_{\tau},1)\prod_i|\psi'(\tau)|^{-\triangle_{\alpha_i}}.
 \end{align*}
Gathering all the above considerations, this proves all our claims.\qed

\section{The case $\gamma=2$}\label{sec:crit}
The case $\gamma=2$ requires some special care because it does not enter the framework of standard Gaussian multiplicative chaos theory. Indeed, recall that by \cite{cf:Kah} we have almost surely
\begin{equation}
\lim_{\epsilon\to 0}\epsilon^2e^{2 X_{\tau,\epsilon}}\,d\lambda_{\hat{g}_\tau}\to 0.
\end{equation}
Formally, we know (see \cite{Rnew7,Rnew12}) that an extra push of order $(\ln\frac{1}{\epsilon})^{1/2}$ is necessary to make this measure converge towards a non trivial limit. Yet, the case of GFF on tori does not exactly enter the framework of \cite{Rnew7,Rnew12}. We explain here how to bridge this gap.

So we consider the family of random measures
\begin{equation}
 M_{2,\tau,\epsilon}:= \sqrt{\pi/2}(\ln\frac{1}{\epsilon})^{1/2}  \epsilon^{2}e^{2 (X_{\tau,\epsilon}-\ln{\rm Im}(\tau)) } \,d\lambda_{\hat{g}_\tau}.
\end{equation}
We claim
\begin{proposition}\label{lawcrit1}
The  following limit exists  in probability
\begin{equation}
M_{2,\tau}:=\lim_{\epsilon\to 0}M_{2,\tau,\epsilon}
=\sqrt{\pi/2}e^{2\Theta_{\tau} - 2\ln{\rm Im}(\tau)} \lim_{\epsilon\to 0}(\ln\frac{1}{\epsilon})^{1/2}  e^{2 X_{\tau,\epsilon}-2\E[ X_{\tau,\epsilon}^2] } \,d\lambda_{\hat{g}_\tau}
\end{equation}
in the sense of weak convergence of measures (with  
$\Theta_{\tau} $ defined in \ref{Thetatau}).
Almost surely, this limiting measure is non trivial, without atom and is   (up to the multiplicative constant $\sqrt{\pi/2}e^{2\Theta_{\tau} -2\ln{\rm Im}(\tau)}$) a critical Gaussian multiplicative chaos of the field $X_{\tau}$ with respect to the measure $\lambda_{\hat{g}_\tau}$. This measure has negative moments of all orders.

Furthermore, for any modular transformation $\psi\in\mathcal{M}$, we have
\begin{equation}
(X_{\psi(\tau)}, M_{2,\psi(\tau)})\stackrel{law}{=}(X_{\tau}\circ\widetilde{\psi }  , \widetilde{\psi}^{-1}*M_{2,\tau}),
\end{equation}
where $ \widetilde{\psi}^{-1}*M_{2,\tau}$ is the pushforward measure of $M_{2,\tau}$ by $ \widetilde{\psi}^{-1}$.
\end{proposition} 

It remains to show that the random measure gives finite mass to $\T$ almost surely (notice that the expectation is infinite)
\begin{proposition}\label{lawcrit2}
Almost surely we have $M_{2, \tau}(\T)<+\infty$.
\end{proposition}

From now on we fix $\mu>0$ and $\gamma=2$, yielding $Q=2$. We give now the definition function of the $\epsilon$-regularized partition function in the metric $g_\tau$ for some fixed $\tau\in\H$  for every measurable bounded function $F$ 
 \begin{align}\label{actiontoricrit}
  \Pi_{2, \mu}^{(z_i,\alpha_i)_i }  ( g_\tau ,F,\epsilon)   := & \\
 Z^{{\rm FF}}( g_\tau) \int_\R  \E &\left[  F\left( c+  X_{g_\tau,\epsilon}  -\frac{Q}{2}\ln{\rm Im}(\tau)+\frac{Q}{2} \varphi_\tau   ,  \sqrt{\pi/2}(\ln\frac{1}{\epsilon})^{1/2} \epsilon^{2}  e^{2 ( c+X_{g_\tau,\epsilon}-\frac{Q}{2}\ln{\rm Im}(\tau)+\frac{Q}{2} \varphi_\tau)}  d\lambda_{\hat{g}_\tau}   \right) \right. \nonumber \\  
 & \quad\prod_i \epsilon^{\frac{\alpha_i^2}{2}}e^{\alpha_i (c+  X_{\tau,\epsilon} (z_i) - \ln{\rm Im}(\tau)+  \varphi^\epsilon_\tau(z_i))} 
 \nonumber\\
 & \left. \exp\Big( -\frac{1}{2\pi}\int_{\T }R_{g_\tau}(c+X_{\tau}  )\,\lambda_{g_\tau} - \mu \sqrt{\pi/2}(\ln\frac{1}{\epsilon})^{1/2} \epsilon^{2}\int_{\T} e^{2 ( c+X_{\tau,\epsilon}- \ln{\rm Im}(\tau)+  \varphi^\epsilon_\tau)}  \,d\lambda_{\hat{g}_\tau}   \Big) \right]\,dc.\nonumber
\end{align}
From this stage onwards, the properties of LQG on tori with $\gamma=2$ as well as their proofs is exactly the same as for $\gamma=2$, namely all the theorems we stated remain valid if we take $\gamma=2$.

\medskip

 \noindent {\it Proof of Proposition \ref{lawcrit1}.} We will construct the random measure in the interior of $\T$ and then show that this measure gives finite mass to $\T$. 
 Observe first that the GFF $X_\tau$ inside $\T$ can be decomposed as 
\begin{equation}
X_\tau=X_{\tau,{\rm Dir}}+P^{\tau}
\end{equation}
where $X_{\tau,{\rm Dir}}$ and $P^{\tau}$ are independent, $X_{\tau,{\rm Dir}}$ is a GFF on $\T$ in the metric $\hat{g}_\tau$ with Dirichlet boundary condition on $\partial \T$ and $P^{\tau}$ is the harmonic extension (in the metric $\hat{g}_\tau$) of the restriction of $X_\tau$ to the boundary of $\T$. In particular, $P^{\tau}$ is a smooth Gaussian process in the interior of $\T$. We denote by $X_{\tau,{\rm Dir},\epsilon}$ and $P^{\tau,\epsilon}$ the $\epsilon$-circle average of the two processes.

The crucial lemma is the following
 \begin{lemma}\label{lawcrit3}
The family of random measures  
$$(\ln\frac{1}{\epsilon})^{1/2} e^{2 X_{\tau,{\rm Dir},\epsilon}-2\E[X_{\tau,{\rm Dir},\epsilon}^2]}\,d\lambda_{\hat{g}_\tau}$$
converges in probability in the sense of weak convergence of measures towards a non trivial limit, which is atom free and possesses negative moments of all order. Furthermore, almost surely
\end{lemma}

Admitting for a while this lemma, it is clear that the random measure $ M_{2,\tau,\epsilon}$ converges in probability in the interior of $\T$ as it can be rewritten as
\begin{align*}
M_{2,\tau,\epsilon}
=&\sqrt{\pi/2}e^{2\E[ X_{\tau,{\rm Dir},\epsilon}^2]+2 \ln\epsilon} \lim_{\epsilon\to 0}(\ln\frac{1}{\epsilon})^{1/2}e^{2P^{\tau,\epsilon}-2\E[(P^{\tau,\epsilon})^2]}  e^{2 X_{\tau,{\rm Dir},\epsilon}-2\E[ X_{\tau,{\rm Dir},\epsilon}^2] } \,d\lambda_{\hat{g}_\tau}.
\end{align*}
Proposition \ref{circlegreentorus} and the fact that the Gaussian process $P^{\tau }$ is continuous inside $\T$ allows us to conclude easily.\qed

\medskip

\noindent {\it Proof of Lemma \ref{lawcrit3} and Proposition \ref{lawcrit2}.} The point is that the convergence is established in \cite{Rnew12} with a white noise decomposition of the GFF $X_{\tau,{\rm Dir}}$ instead of the circle average regularization we are dealing with. We explain now how the convergence for white noise decomposition implies convergence of the circle average approximations.

Let us denote by $\bar{X}_{\tau}$ another GFF with Dirichlet boundary condition on $\partial \T$, independent of everything,  admitting a white noise decomposition based on the heat kernel of the Brownian motion killed upon touching the boundary $\partial \T$ constructed in the manner explained in \cite[section 6]{Rnew12}. Let us denote by   $(\bar{X}^{wn}_{\tau,\epsilon})_{\epsilon}$ the corresponding white noise approximations, which are measurable with respect to the underlying white noise $W$ distributed on $[1,+\infty[\times\T$. From  \cite[section 6]{Rnew12}, the family of random measures $$\bar{M}_\epsilon:=(\ln\frac{1}{\epsilon})^{1/2} e^{2 \bar{X}^{wn}_{\tau,\epsilon}-2\E[(\bar{X}^{wn}_{\tau,\epsilon})^2]}\,d\lambda_{\hat{g}_\tau}$$ converges in probability in the sense of weak convergence on measures on $\T$ towards a limiting random measures $\bar{M}$ possessing all the desired properties, which is a measurable function of the white noise $W$, i.e. $\bar{M}=F(W)$ for some measurable function $F$. Furthermore from \cite[section 6]{Rnew12}, we have almost surely $\bar{M}(\T)<+\infty$. Hence, denoting by $\bar{M}'$ the measure on $\T$ defined by
$$\bar{M}' =e^{2P^{\tau}-2\E[(P^{\tau})^2]}\,d\bar{M}$$ we have 
$$\E[\bar{M}'(\T)|\bar{M}]=\bar{M}(\T)<+\infty.$$
As we will show below that the measure $M_{2,\tau}$ has the same law as $\bar{M}'$ up to a deterministic multiplicative factor, this will show Proposition \ref{lawcrit2}.

Let us denote by $(\bar{X}^{ca}_{\tau,\epsilon})_{\epsilon}$ the circle average regularization of the GFF $\bar{X}_{\tau}$. Now, we show convergence in probability of $(\ln \frac{1}{\epsilon})^{1/2}\, e^{ 2\bar{X}^{ca}_{\tau,\epsilon} -2\E[(\bar{X}^{ca}_{\tau,\epsilon})^2]}  d\lambda_{\hat{g}_\tau}$  towards the same limit $\bar{M}$. The ideas in the following stem from the techniques developed in \cite{review} along with some some variant of lemma 49 in \cite{shamov} (we will not recall lemma 49 as our proof will be self contained). For this, we introduce an independent copy $\widetilde{X}_{\tau}$ of $\bar{X}_{\tau}$, and $(\widetilde{X}^{ca}_{\tau,\epsilon})_{\epsilon}$ its circle average approximation.  Let us define for $t\in[0,1]$ and $x\in \T$
$$
Z_\epsilon(t,x)=\sqrt{t}\widetilde{X}^{ca}_{\tau,\epsilon}(x)+\sqrt{1-t} \bar{X}^{ca}_{\tau,\epsilon}(x).
$$       
Now, we set 
\begin{equation*}
\widetilde{M}_{ \epsilon}=(\ln \frac{1}{\epsilon})^{1/2}\, e^{ 2\widetilde{X}^{ca}_{\tau,\epsilon} -2\E[(\widetilde{X}^{ca}_{\tau,\epsilon})^2]}  d\lambda_{\hat{g}_\tau}.
\end{equation*}

We first show that the random measures $(\widetilde{M}_{ \epsilon})_{\epsilon}$ converge in law to $\bar{M}=F(W)$. From \cite[Proof of Theorem 2.1]{vincent}, one gets that for all $0<\alpha<1$ and ball $B\subset \T$ such that ${\rm dis}(B,\partial \T)>0$
\begin{align*}
& \underset{\epsilon \to 0} {\overline{\lim}}  \big|\E[\widetilde{M}_{ \epsilon}(B)^\alpha] -\E[\bar{M}_{ \epsilon}(B)^\alpha]\big|  \\
& \leq c \frac{\alpha(1-\alpha)}{2}C_A \underset{\epsilon \to 0} {\overline{\lim}} \int_0^1\E\Big[\Big ( (\ln \frac{1}{\epsilon})^{1/2} \int_{\T} e^{2Z_\epsilon(t,\cdot)-2\E[Z_\epsilon(t,\cdot)^2]}\,d\lambda_{\partial}  \Big )^\alpha\Big] \,dt \\
& +  c\,  \overline{C}_A  \underset{\epsilon \to 0} {\overline{\lim}} \int_0^1 \E \Big [\Big ( \sup_{0\leq i<\frac{1}{A\epsilon}}  (\ln \frac{1}{\epsilon})^{1/2}  \int_{[2iA\epsilon,2(i+1)A\epsilon]^2\cap B}e^{2Z_\epsilon(t,x)-2\E[Z_\epsilon(t,x)^2]}  d\theta \Big)^\alpha     \Big]  dt,   
\end{align*}
where
\begin{equation*}
C_A=\underset{\epsilon \to 0}{\overline{\lim}} \sup_{|x-x'| \geq A \epsilon}  | \E[ \widetilde{X}^{ca}_{\tau,\epsilon}(x)\widetilde{X}^{ca}_{\tau,\epsilon}(x')   ] -  \E[  \bar{X}^{wn}_{\tau,\epsilon}(x)  \bar{X}^{wn}_{\tau,\epsilon}(x')   ] |
\end{equation*}
and
\begin{equation*}
\overline{C}_A= \underset{\epsilon \to 0}{\overline{\lim}} \sup_{|x-x'| \leq A \epsilon}  | \E[ \widetilde{X}^{ca}_{\tau,\epsilon}(x) \widetilde{X}^{ca}_{\tau,\epsilon}(x')    ] -  \E[\bar{X}^{wn}_{\tau,\epsilon}(x) \bar{X}^{wn}_{\tau,\epsilon}(x')    ] |.
\end{equation*}
The reader can check that $\overline{C}_A$ is bounded independently of $A$ and $\underset{A \to \infty}{\lim} C_A=0$. Since \\ $\E\big[ \left ((\ln \frac{1}{\epsilon})^{1/2} \int_\T e^{2Z_\epsilon(t,x)-2\E[Z_\epsilon(t,x)^2]}\,\lambda_{\hat{g}_\tau} (dx) \right )^\alpha\big] $ is also bounded  independently of everything for any fixed $\tau$ (by comparison with Mandelbrot's multiplicative cascades as explained in the \cite[appendix]{Rnew7} and \cite[appendix B.4]{Rnew12}), we are done if we can show that for all $t\in [0,1]$
\begin{equation}\label{claim}
 \underset{\epsilon \to 0} {\overline{\lim}} \,\,\E \Big [\Big ( \sup_{0\leq i<\frac{1}{A\epsilon}} (\ln \frac{1}{\epsilon})^{1/2} \int_{[2iA\epsilon,2(i+1)A\epsilon]^2\cap B}e^{2Z_\epsilon(t,x)-2\E[Z_\epsilon(t,x)^2]}  \,\lambda_{\hat{g}_\tau} (dx) \Big)^\alpha     \Big] =0.
 \end{equation}
Notice that this quantity is less than
\begin{equation}\label{ascola}
\big(\ln \frac{1}{\epsilon}\big)^{\alpha/2}\epsilon^\alpha \E  \Big [  \Big  (    e^{\sup_{x\in B} 2Z_\epsilon(t,x)-2\E[Z_\epsilon(t,x)^2]}   \Big )^\alpha     \Big ] .
\end{equation} 
To estimate this quantity, we use the main result  of \cite{acosta}: more precisely, setting 
$$m_\epsilon=2\ln\frac{1}{\epsilon}-\frac{3}{2}\ln\ln \frac{1}{\epsilon},$$
 we claim that there exist two constants $C,c>0$ such that for $\epsilon$ small enough
\begin{align*}
\forall v\geq 0,\quad \P\Big(\Big|\max_{x\in B}2Z_\epsilon(t,x)-m_\epsilon \Big|\geq v\Big)\leq Ce^{-cv}.
\end{align*}
In particular we get that for   $\alpha<c$  
\begin{equation*}
\sup_{\epsilon} \E  \Big [  \Big  (    e^{\sup_{x\in B} 2Z_\epsilon(t,x)}   \Big )^\alpha     \Big ]<\infty.
\end{equation*}
Plugging this estimate into \eqref{ascola}, we see that the quantity \eqref{ascola} is less than
$$C' \big(\ln \frac{1}{\epsilon}\big)^{\alpha/2}\epsilon^{2\alpha}e^{\alpha m_\epsilon}=C' \big(\ln \frac{1}{\epsilon}\big)^{-\alpha}. $$
for some constant $C'>0$. This proves the claim \eqref{claim}, hence the convergence in law of the random measures $(\widetilde{M}_{ \epsilon})_{\epsilon}$ towards  $\bar{M}=F(W)$. 

Now we deduce that the family $(W,\widetilde{M}_{ \epsilon})_\epsilon$ converges in law. Take any  smooth function $R$ on $[1,+\infty[\times  \T$ with compact support in the interior of $\T$, any  continuous function $g$ with compact support in the interior of $\T$ and any bounded continuous function $G$ on $\R$. We have by using the Girsanov transform 
\begin{equation}\label{WR}
\E[e^{W(R)}G(\widetilde{M}_{ \epsilon}(g))]=e^{\frac{1}{2}{\rm Var}[W(R)]}\E[ G(\widetilde{M}_{ \epsilon}(e^{T_\epsilon(R)}g)]
\end{equation}
where $T_\epsilon(R)$ is defined by $$x\in\T\mapsto T_\epsilon(R)(x):=\E[ \widetilde{X}^{ca}_{\tau,\epsilon}(x)W(R)].$$ This is a  continuous function which converges uniformly as $\epsilon\to 0$ towards the continuous
$$z\in\partial \D\mapsto T(R)(x):=\E[\widetilde{X}_{\tau}(x)W(R)].$$
The quantity in the right-hand side of \eqref{WR} converges as $\epsilon\to 0$ towards
 $$e^{\frac{1}{2}{\rm Var}[W(R)]}\E[ G(\bar{M}(e^{T(R)}g)]=\E[e^{W(R)}G(\bar{M}(g))].$$
Hence our claim about the convergence in law of the couple $(W,\widetilde{M}_{ \epsilon})_\epsilon$ towards $(W, \bar{M}=F(W))$.
 
 Now we consider the family $(W, \widetilde{M}_{ \epsilon},F(W))_\epsilon$, which is tight. Even if it means extracting a subsequence, it converges in law towards some $(\mathcal{W},  \mathcal{M}, \bar{\mathcal{M}})$. We have just shown that the law of $(\mathcal{W},  \mathcal{M})$ is that of  $(\mathcal{W},  F
(\mathcal{W}))$, i.e. the same as the law of 
$(\mathcal{W}, \bar{\mathcal{M}})$. Hence $\mathcal{M} =\bar{\mathcal{M}}$ almost surely. Therefore $\widetilde{M}_{ \epsilon}- F(W)$ converges in law towards $0$, hence in probability.   Since the  convergence in probability of the family $(\widetilde{M}_{ \epsilon})_\epsilon$ implies the convergence of probability of every family  that has the same law as $(\widetilde{M}_{ \epsilon})_\epsilon$, this entails the convergence in probability of the random measures $(\ln\frac{1}{\epsilon})^{1/2} e^{2 X_{\tau,{\rm Dir},\epsilon}-2\E[X_{\tau,{\rm Dir},\epsilon}^2]}\,d\lambda_{\hat{g}_\tau}$ towards some limit that has the same law as $\bar{M}$. \qed

\section{LQG on tori: law of the Liouville modulus and conjecture related to random planar maps}\label{sec:LQG}

\subsection{Reminders on LQFT and LQG}
\label{LQFT&LQG}
Before stating our result, let us briefly (and softly) recall to mathematicians how one builds Liouville quantum gravity (LQG) out of the Liouville QFT (on manifolds with fixed topology, here the torus). It would be more convenient for the reader to be familiar with the basic axiomatic of conformal field theories (CFT, see for instance \cite[Section 2]{gaw}) before reading this section. 

The LQG partition function is the product   of three independent conformal field theory (CFT): matter+ghost+LQFT (see \cite{Pol,friedan}). We will not explain in further details what these objects correspond to: what will be important for us is just the fact that the partition function of the ghost CFT and (in many cases) the matter field CFT are explicitly known.  These partition functions depend on the background metric on which they are computed. 
 
So, to construct Liouville quantum gravity on tori, we basically need the three ingredients below:
\begin{enumerate}
\item the partition function of the ghost CFT: the ghost CFT is defined for instance in \cite{friedan}, the partition function on the tori is calculated in \cite{polchinski,gupta}
\begin{equation}
\label{partghost}
Z_{{\rm Ghost}}(e^{\varphi_\tau}\hat{g}_\tau)
= e^{-\frac{26}{96\pi}\int_\T|\partial^\tau\varphi_\tau|_\tau^2\,d\lambda_\tau}\frac{|\eta(\tau)|^4}{2\,{\rm Im}(\tau)}.
\end{equation}
\item the partition function of LQFT: it is nothing but the quantity $\Pi^{(z_i,\alpha_i)_i}_{\gamma,\mu}( e^{\varphi_\tau}\hat{g}_\tau)$ studied in this paper.
\item a CFT for the matter field: for each metric $g_\tau $  on $\T$, it consists in the partition function $Z_{{\rm Matter}}(g_\tau)$ and a set of correlation functions   for the``primary"  matter fields $\theta_i$ denoted by
$$Z_{{\rm Matter}}(g_\tau,\theta_{i_1}(z_1),\dots,\theta_{i_n}(z_n))$$ where $n\geq 1$ and $i_1,\dots,i_n$ belong to a fixed set  $I$ and $z_1,\dots,z_n$ are distinct points in $\T$. These correlation functions are supposed to obey the following rules
\begin{description}
\item[$\bullet$] (modular covariance) for every $\psi\in\mathcal{M}$, 
$Z_{{\rm Matter}}(\hat{g}_{\psi(\tau)})=Z_{{\rm Matter}}(\hat{g}_\tau)$ and 
\begin{equation}\label{mattermodcov}
 \frac{Z_{{\rm Matter}}(\hat{g}_{\psi(\tau)},\theta_{i_1}(z_1),\dots,\theta_{i_n}(z_n))}{Z_{{\rm Matter}}(g_\tau,\theta_{i_1}(\widetilde{\psi}(z_1)),\dots,\theta_{i_n}(\widetilde{\psi}(z_n)))}=\Big(\prod_{k=1}^n|\psi'(\tau)|^{-\triangle^{{\rm m}}_{i_k}}\Big),
\end{equation}
 where $\triangle^{{\rm m}}_{i_k}$ is the conformal weight  of the primary matter field $\theta_{i_k}$ (we restrict ourselves to the case of rational CFT and to spin $s=0$ primary fields such that $\triangle^{{\rm m}}_{i_k}$ is real).
\item[$\bullet$] (Weyl anomaly) for every metric $e^{\varphi_\tau}\hat{g}_\tau$ conformally equivalent to $\hat{g}_\tau$
\begin{equation}\label{partmatter}
Z_{{\rm Matter}}(e^{\varphi}\hat{g}_\tau)=e^{\frac{c_{{\rm m}}}{96\pi}\int_\T|\partial^\tau\varphi_\tau|_\tau^2  \,d\lambda_\tau}Z_{{\rm Matter}}( \hat{g}_\tau)
\end{equation}
and 
\begin{equation}\label{corrmatter}
\frac{Z_{{\rm Matter}}(e^{\varphi}\hat{g}_\tau ,\theta_{i_1}(z_1),\dots,\theta_{i_n}(z_n))}{Z_{{\rm Matter}}( \hat{g}_\tau,\theta_{i_1}(z_1),\dots,\theta_{i_n}(z_n))}=e^{\frac{c_{{\rm m}}}{96\pi}\int_\T|\partial^\tau\varphi_\tau|_\tau^2  \,d\lambda_\tau} 
\end{equation}
 where $c_{{\rm m}}\leq 1$ is the central charge of the matter field CFT.
\end{description}
\end{enumerate}

\begin{remark}
The most elementary example of matter field for CFT that the reader can have in mind is that of the compactified boson explained in \cite{gaw}.  
\end{remark}  

\begin{remark}{{\bf (Conformal ansatz)}}
Let us explain in passing an important concept in physics called the conformal ansatz \cite{Pol,cf:Da,DistKa}, which fixes the value of $Q$ (hence $\gamma$) in terms of $c_{{\rm m}}$. The choice of the metric $\hat{g}_\tau$ is somewhat arbitrary and result from the need of fixing the conformal gauge: i.e we have picked up arbitrarily  a family of representatives $(\hat{g}_\tau)_{\tau}$ of all the conformal structures on $\T$. Yet physics should not depend on this choice and should remain the same had we chosen another family of representatives, call it $(e^{\varphi_\tau}\hat{g}_\tau)_\tau$,  for some modular log-conformal factor $\varphi_\tau$. This means that the partition function for LQG with $n$ marked points $(z_i)_{i}$ is
\begin{equation}
Z_{{\rm Matter}}(e^{\varphi_\tau}\hat{g}_\tau,z_1,\dots,z_n)Z_{{\rm Ghost}}(e^{\varphi_\tau}\hat{g}_\tau)\Pi^{(z_i,\alpha_i)_i}_{\gamma,\mu}(e^{\varphi_\tau}\hat{g}_\tau)
\end{equation}
and should not depend on the conformal factor $e^{\varphi_\tau}$ \cite{Pol,cf:KPZ,cf:Da,DistKa}. By using the  Weyl anomalies of Theorem \ref{th:weyl} and \eqref{partghost}$+$\eqref{corrmatter}, we see that this product is equal to 
\begin{equation}
e^{\frac{c_{{\rm m}}+1+6Q^2-26}{96\pi}\int_\T|\partial^\tau\varphi_\tau|_\tau^2\,d\lambda_\tau} Z_{{\rm Matter}}( \hat{g}_\tau,z_1,\dots,z_n)Z_{{\rm Ghost}}( \hat{g}_\tau)\Pi^{(z_i,\alpha_i)_i}_{\gamma,\mu}( \hat{g}_\tau).
\end{equation}
Therefore the quantity $c_{{\rm m}}+1+6Q^2-26$ must vanish and this directly gives the famous KPZ relation \cite{cf:KPZ}
\begin{equation}\label{string}
\gamma=\frac{\sqrt{25-c_{{\rm m}}}-\sqrt{1-c_{{\rm m}}}}{\sqrt{6}} \quad \text{ with  }\gamma \in ]0,2]\text{ for }c_{{\rm m}}\leq 1.
\end{equation}
Furthermore, for modular invariance to be preserved, the $(\alpha_i)_i$ (satisfying the Seiberg bounds) must be chosen in such a way that  $\triangle^{{\rm m}}_{i_k}+\triangle_{\alpha_i}=1$, that is 
\begin{equation}\label{matterKPZ}
\triangle^{{\rm m}}_{i_k}+\frac{\alpha}{2}(Q-\frac{\alpha}{2})=1.
\end{equation}
When $\triangle^{{\rm m}}_{i_k}=0$, this means that $\alpha_i=\gamma$.
\end{remark}

\subsection{Law of the Liouville modulus}
\label{Liouvillemodulus}
Let us go now back to the main point of this section. In what follows, we will assume that we are given the $n$-point correlation function of a matter field CFT $Z_{{\rm Matter}}$ with $c_{{\rm m}}\leq 1$ (i.e. the content of 2 above). We fix the value of $\gamma\in]0,2]$ according to the relation \eqref{string} in such a way that we may assume that the conformal factor  $ \varphi_\tau$ vanish, i.e. $\varphi_\tau=0$,  because of Weyl invariance.

We define the $n$-correlation function of LQG   on tori with parameters ${\rm c}_{{\rm m}},\mu$ and $n$ primary fields $(\theta_{i_k})_k$ at the points $z_1,\dots,z_n$  applied to a bounded continuous functional $G$ on $H^{-1}(\T)\times M_+(\T)\times \mathcal{S}$ by
$$Z^{ (z_k,\theta_{i_k})_k}_{{\rm c}_{{\rm m}},\mu,{\rm LQG}}(G)=\int_{\mathcal{S}} Z_{{\rm Matter}}(\hat{g}_{\tau},\theta_{i_1}(z_1),\dots,\theta_{i_n}(z_n))\Pi_{\gamma,\mu}^{(z_k,\alpha_k)_k}(\hat{g}_\tau,G(.,\tau)) {\rm Im}(\tau)^n Z_{{\rm Ghost}}(\hat{g}_\tau)d\tau .$$
From Theorem \ref{th:seibergtau} as well as \eqref{partGFF}+\eqref{partghost}, it can be rewritten as 
\begin{align}
Z^{ (z_k,\theta_{i_k})_k}_{{\rm c}_m,\mu,{\rm LQG}} (G)\label{LQGexp}&\ = \\
&\frac{1}{2}\int_{\mathcal{S}}  \int_\R  e^{c\sum_k\alpha_k}\E\Big[G\big(X_\tau-\frac{Q}{2}\ln{\rm Im}(\tau)+H_\tau+c,e^{\gamma c+\gamma H_\tau}\,dM_{\gamma,\tau},\tau\big)e^{-\mu e^{\gamma c}\int_\T e^{\gamma H_\tau}\,dM_{\gamma,\tau}}\Big]\nonumber\\
&\times   e^{C_\tau(\mathbf{z})}Z_{{\rm Matter}}(\hat{g}_{\tau},\theta_{i_1}(z_1),\dots,\theta_{i_n}(z_n)){\rm Im}(\tau)^n  \sqrt{{\rm Im}(\tau)} |\eta(\tau)|^2    \lambda_{\mathcal{S}}(d\tau)  \ dc\nonumber.
\end{align}
provided that this quantity is finite. In what follows, we wish to see  the modulus $\tau$ as a random variable. Therefore we have to deal with positive quantities. So we assume in the following that the matter field correlation functions are positive and that the above integral is finite when $G=1$ (this is for instance the case  for the correlation functions of the compactified boson with imaginary background charge \cite{gaw}).  Our purpose below is to find an explicit expression for the probability density of the Liouville modulus.
\begin{definition}{\bf (Law of the Liouville measure/field/modulus)}
The joint law of the Liouville field, measure and modulus of LQG on tori with parameters ${\rm c}_m,\mu$ and $n$ primary fields $(\theta_{i_k})_k$ at the points $z_1,\dots,z_n$  is defined by 
\begin{equation}
E^{(z_k,\theta_{i_k})_k}_{{\rm c}_m,\mu,{\rm LQG}}[G(\phi,Z,\tau )]=\frac{Z^{(z_k,\theta_{i_k})_k}_{{\rm c}_m,\mu,{\rm LQG}}(G)}{Z^{(z_k,\theta_{i_k})_k}_{{\rm c}_m,\mu,{\rm LQG}}(1)}
\end{equation}
for every bounded continuous functional $G$ on $H^{-1}(\T)\times M_+(\T)\times \mathcal{S}$.
We denote by $\P^{(z_k,\theta_{i_k})_k}_{{\rm c}_m,\mu,{\rm LQG}}$ the corresponding probability measure.
\end{definition}

\begin{proposition}{\bf (Joint law of the Liouville volume/modulus)}\label{joint}
Under $\P^{(z_k,\theta_{i_k})_k}_{{\rm c}_m,\mu,{\rm LQG}}$,  the volume of space $Z(\T)$ and Liouville modulus  $ \tau$ are independent. The volume follows the Gamma law $\Gamma (\gamma^{-1}\sum_i\alpha_i,\mu)$ on $\R_+$, i.e. has density
\begin{equation}
\frac{\mu^{\frac{\sum_i\alpha_i}{\gamma}}}{\Gamma\Big(\frac{\sum_i\alpha_i}{\gamma}\Big)}\ind_{\R_+}(y) y^{\frac{\sum_i\alpha_i}{\gamma}-1}e^{-\mu y}\,dy,
\end{equation}
whereas the Liouville modulus has the probability law   on $\mathcal{S}$
\begin{equation}
\frac{1}{R}  \E\Big[\Big(\int_\T e^{\gamma H_\tau}\,dM_{\gamma,\tau}\Big)^{-\frac{\sum_i\alpha_i}{\gamma}}\Big]e^{C_\tau(\mathbf{z})} Z_{{\rm Matter}}(\hat{g}_{\tau},\theta_{i_1}(z_1),\dots,\theta_{i_n}(z_n))\,{\rm Im}(\tau)^n  \sqrt{{\rm Im}(\tau)} |\eta(\tau)|^2 \lambda_{\mathcal{S}}(d\tau)
\end{equation}
where $R$ is a deterministic renormalizing constant fixed so as to have a probability measure on $\mathcal{S}$. Conditionally on the Liouville modulus $\tau$, the fixed volume LQG measure has law
\begin{equation}
\E^{n,(\theta_{i_k})_k}_{{\rm c}_m,\mu,{\rm LQG}}[G(Z(dx))|Z(\T)=y,\tau]=\frac{ \E\Big[G\big( y\frac{e^{\gamma H_\tau}\,dM_{\gamma,\tau}}{\int_\T e^{\gamma H_\tau}\,dM_{\gamma,\tau}} \big) \Big(\int_\T e^{\gamma H_\tau}\,dM_{\gamma,\tau}\Big)^{-\frac{\sum_i\alpha_i}{\gamma}}\Big] }{ \E\Big[ \Big(\int_\T e^{\gamma H_\tau}\,dM_{\gamma,\tau}\Big)^{-\frac{\sum_i\alpha_i}{\gamma}}\Big]  }.
\end{equation}
\end{proposition}

 \noindent {\it Proof.} Let us start from the relation \eqref{LQGexp} and make the change of variables $y=e^{\gamma c} \int_\T e^{\gamma H_\tau}\,dM_{\gamma,\tau}$ to get (for some $G$ which does not depend on the Liouville field)
\begin{align}
Z^{(z_k,\theta_{i_k})_k}_{{\rm c}_m,\mu,{\rm LQG}}(G)=&\frac{1}{2\gamma }\int_{\mathcal{S}}  \int_{\R_+} y^{\frac{1}{\gamma}\big(\sum_i\alpha_i\big)-1}  \E\Big[G\big( y\frac{e^{\gamma H_\tau}\,dM_{\gamma,\tau}}{\int_\T e^{\gamma H_\tau}\,dM_{\gamma,\tau}},\tau\big) \Big(\int_\T e^{\gamma H_\tau}\,dM_{\gamma,\tau}\Big)^{-\frac{\sum_i\alpha_i}{\gamma}}\Big]e^{-\mu y} \nonumber\\
&\times   e^{C_\tau(\mathbf{z})}Z_{{\rm Matter}}(\hat{g}_{\tau},\theta_{i_1}(z_1),\dots,\theta_{i_n}(z_n)) {\rm Im}(\tau)^n  \sqrt{{\rm Im}(\tau)} |\eta(\tau)|^2     \lambda_{\mathcal{S}}(d\tau)dy.\label{LQG2}
\end{align}
From this relation we get that
\begin{align*}
\E^{(z_k,\theta_{i_k})_k}_{{\rm c}_m,\mu,{\rm LQG}}&[G(Z(\T),\tau )]\\
=& \frac{1}{Z^{(z_i,\alpha_i)_i}_{{\rm c}_m,\mu,{\rm LQG}}(1)}\int_{\mathcal{S}}  \int_{\R_+} G\big( y,\tau\big) y^{\frac{1}{\gamma}\big(\sum_i\alpha_i\big)-1}\Big(\int_{\T^n} \E\Big[ \Big(\int_\T e^{\gamma H_\tau}\,dM_{\gamma,\tau}\Big)^{-\frac{\sum_i\alpha_i}{\gamma}}\Big]e^{-\mu y} \\
&\times   e^{C_\tau(\mathbf{z})}Z_{{\rm Matter}}(\hat{g}_{\tau},\theta_{i_1}(z_1),\dots,\theta_{i_n}(z_n)) {\rm Im}(\tau)^n  \sqrt{{\rm Im}(\tau)} |\eta(\tau)|^2     \lambda_{\mathcal{S}}(d\tau).
\end{align*}
This proves our first claim about the independence of the volume and the modulus together with their respective laws. The second claim about the conditional law is then a direct consequence of \eqref{LQG2}.\qed


\subsection{Conjectures related to genus$=1$ rooted 
maps and quadrangulations
}
\label{sec:conj}

In this section we will  present some precise conjectures on the connection of our results to the work on discrete models of 2d gravity, randoms surfaces and random planar maps. The standard way to discretize 2d quantum gravity coupled to matter fields is to consider a statistical mechanics model (corresponding ai its critical point to a conformal field theory with central charge $c_\textrm{m}$), defined on a  random lattice, corresponding to the random metric.

We formulate below precise mathematical conjectures on the relationship of LQG to that setup in the simplest case of pure gravity, i.e. when there is no coupling with matter. The reader may easily extend our conjectures to the general case by adapting to the torus the general picture drawn in \cite{DKRV}.

There are numerous results on random lattices with the topology of the torus, from the first results of \cite{bessis} using the random matrices techniques pioneered in \cite{BIPZ} to recent results of for instance \cite{bett1,bett2,bender, Chapuy, miermont}. For simplicity we shall discuss only quadrangulations but these considerations can be straightforwardly extended to more general maps (bipartite quadrangulations, triangulations, etc.) .

Let $\mathcal{Q}_{N}^1$ be the set of rooted quadrangulations  with $N$ faces with the topology of the torus. From  \cite{bender}, we have
$$\mathcal{Q}_{N}^1\sim C\,A^N  $$
 for some constants $A,\,C>0$. For bipartite triangulations $A=12$.
 
Now to each quadrangulation $Q\in \mathcal{Q}_{N}^1$   we associate a standard conformal structure (by gluing Euclidean squares along their edges as prescribed by the quadrangulation). Such a complex manifold is conformally equivalent  to $\T$ equipped with the metric $\hat{g}_\tau$ for some unique $\tau\in\mathcal{S}$. So we consider a conformal map sending this quadrangulation equipped with its complex structure to $(\T,\hat{g}_\tau)$  such that the root gets mapped to $0$. We give volume $a^2$ to each quadrilateral and we denote $\nu_{Q,a}$ the corresponding volume measure on $\T$.

\medskip   
For $\bar{\mu}>\bar{\mu}_c=\ln A$, the full partition of the rooted quadrangulations of the torus reads
 \begin{equation}
\label{Zmub}
Z_{\bar{\mu}}=\sum_Ne^{-\bar{\mu} N}|\mathcal{Q}_{N}^1|
\end{equation}
 converges and we can sample a random quadrangulation according to this partition function. We are interested in the regime where the system samples preferably the quadrangulations with a large number of faces.    Therefore,  we are interested in the limit $\bar{\mu}\to \bar{\mu}_c$ in the following regime: we assume that $\bar{\mu}$ depends on a parameter   $a>0$ such that 
\begin{equation}\label{cosmoRPM}
\overline{\mu}=\bar{\mu}_c+\mu a^2
\end{equation}
 where $\mu$ is a fixed positive constant.  Now, we consider the random measure $\nu_{a,\overline{\mu}}$ on the tori defined by 
\begin{equation*}
\E^{a,\overline{\mu}}[  F( \nu_{a,\overline{\mu}} )  ]= \frac{1}{Z_{a}}\sum_N e^{-(\overline{\mu} - \bar{\mu}_c) N}\sum_{Q \in \mathcal{Q}_{N}^1} F(\nu_{T,a} ),
\end{equation*} 
 for positive bounded functions $F$ where $Z_a$ is a normalization constant. We denote by $\P^{a,\overline{\mu}}$ the probability law associated to $\E^{a,\overline{\mu}}$.

\medskip
We can now state a precise mathematical conjecture:

\begin{conjecture}\label{conjcartes}
Under $\P^{a,\overline{\mu}}$ and under the relation \eqref{cosmoRPM}, the family of random measures  $(\nu_{a,\overline{\mu}})_{a>0} $ converges in law as $a\to0$ in the space of Radon measures equipped with the topology of weak convergence towards the law of the Liouville measure $Z$     under $\P^{(0,0)}_{c_{{\rm m}}=0,\mu,{\rm LQG}}$ (see next subsection for an explicit description). 
\end{conjecture}
The superscript $(0,0)$ in $\P^{(0,0)}_{c_{{\rm m}}=0,\mu,{\rm LQG}}$ means that we consider the $1$ point correlation function of LQG with $z_1=0$ and a primary matter field that is the identity $\theta_{i_1}=\mathbf{1}$ with dimension $\triangle_{i_1}=0$.

Note that  $\nu_{a,\overline{\mu}}(\S^2 ) $ converges in law under $\P^{a,\overline{\mu}}$ as $a\to 0$ towards a $\Gamma(1,\mu)$ distribution, which corresponds precisely to the law of the volume of the space for  LQG with these parameters (see Proposition \ref{joint} or below). The reader may  consult \cite{ambjorn} for further discussions and numerical simulations on this topic.

\subsection{Explicit expression of the law of subsection \ref{sec:conj}}
We give here a description for all values of the central charge ${\rm c}_{\rm m}<1$ so that we can discuss below the law when random planar maps are weighted by a discretized CFT. Recall the relation \eqref{string} between the central charge and $\gamma$. Notice that $n=1$ and $\theta_{i_1}=\mathbf{1}$ (no primary matter field) so that $\triangle^{{\rm m}}_{i_1}=0$. Formula \eqref{matterKPZ} then gives $\alpha_1=\gamma$. Since $\theta_{i_1}=\mathbf{1}$, we only need to know the partition function for the CFT on the torus, $Z_{{\rm Matter}}(\hat{g}_\tau)$  (that we denote $Z_{{\rm Matter}}(\tau)$ when there is no confusion).
This partition function depends explicitly on the operator content of the theory, and of the fusion rules between the primary operators, see \cite{difrancesco} for a detailed discussion. 
In the special case (of interests for string theory) where the matter fields are $D$ copies of GFF, the central charge of the matter sector, and the partition function are simply
\cite{friedan,gupta}
\begin{equation}\label{pm}
c_{{\rm m}}=D  \quad,\qquad Z_{{\rm Matter}}(\hat{g}_\tau)=  {Z^{{\rm FF}}(\tau)}^{c_{{\rm m}}}   \ , 
\end{equation}%
 for some positive constant $C$. 
 In the general case $Z_{{\rm Matter}}(\hat{g}_\tau)$ is more complicated, but it is always a modular invariant function (under transformations of $\tau$ by the modular group $\mathcal{M}$), and with the correct transformation law under conformal changes of the metric $\hat{g}_\tau$ discussed in \ref{LQFT&LQG}.

Using the definition of the law of the Liouville measure and the formula \eqref{pm}, one has
\begin{align}\label{lawmodulus1}
 \E^{(0,0)}_{{\rm c}_{\rm m},\mu,{\rm LQG}}[G(Z,\tau)]
 =& \frac{1}{R}\int_{\R_+} \int_{\mathcal{S}}  \E\Big[ G\Big(y\frac{  e^{\gamma H_\tau}\,dM_{\gamma,\tau}}{\int_\T e^{\gamma H_\tau }\,dM_{\gamma,\tau}},\tau\Big)\Big(\int_\T e^{\gamma H_\tau }\,dM_{\gamma,\tau}\Big)^{-1}\Big] \\
 \quad & {\rm Im}(\tau)\  e^{\frac{\Theta_\tau}{2}\gamma^2-\frac{Q\gamma}{2}\ln{\rm Im}(\tau)}   Z_{{\rm Matter}}(\tau)   \sqrt{{\rm Im}(\tau)} |\eta(\tau)|^2  \lambda_{\mathcal{S}}(d\tau)e^{-\mu y}\,dy,\nonumber
\end{align}

where $R$ is a renormalization constant so as to deal with a probability measure. and $H_\tau(z)=\gamma G_\tau(z )$. 
Observe that this quantity is modular invariant. It is then readily seen that
 \begin{itemize}
\item the law of the volume of space $Z(\T)$ follows the law $\Gamma(1,\mu)$.
\item the law of the Liouville modulus is given by (up to constant to deal with probability law) $$Z_{{\rm Matter}}(\tau)\,  \sqrt{{\rm Im}(\tau)} |\eta(\tau)|^2  \lambda_{\mathcal{S}}(d\tau)$$  
thereby recovering the formula in \cite{ambjorn,gupta} in the special case \ref{pm}.
\end{itemize}

\medskip We precise here a few examples of different CFT with various central charge 
\begin{itemize}
\item {\bf Pure gravity}: this is the case when the CFT has central charge $c_{{\rm m}}=0$, for instance when there is no matter field coupled to gravity (in which case $Z_{{\rm Matter}}(\tau)=1$) or the scaling limit of $2d$ critical percolation coupled to gravity. Then \eqref{string} shows that $\gamma=\sqrt{8/3}$.  
\item {\bf Ising model}: the Ising model at criticality has central charge $c_{{\rm m}}=\frac{1}{2}$, in which case \eqref{string} gives $\gamma=\sqrt{3}$. In that case, $Z_{{\rm Matter}}(\tau)=\big|\frac{\vartheta_2(0,\tau)}{2 \eta(\tau)}\big|+\big|\frac{\vartheta_3(0,\tau)}{2 \eta(\tau)}\big|+\big|\frac{\vartheta_4(0,\tau)}{2 \eta(\tau)}\big|$, where $\vartheta_2$, $\vartheta_3$, $\vartheta_4$ are the auxiliary theta functions (see \cite{difrancesco}).
\item {\bf Compactified boson with imaginary background charge $E$}: it can be seen as a suitable shifted GFF with values in $\R\setminus (2\pi\Z)$  (see \cite[lecture 1 section 4]{gaw} or \cite[section 2.1.3]{dubedattorsion}). The central charge is $c_{{\rm m}}=1-6E^2$ for $E\geq 0$, in which case \eqref{string} gives $\gamma=\sqrt{4+E^2}-E$.
\end{itemize}
Notice however that for ${\rm c}_{\rm m}=1$, the vertex operator $e^{2 X_\tau}$ does not satisfy the Seiberg bound so that one must adapt our conjectures in the case of random planar maps weighted by CFT with central charge $1$.

 \appendix
   
\section{Special functions}\label{special}
The reader may consult \cite{whit} for the content of this section. We consider the special functions: 

\noindent $\bullet$ The {\it Dedekind's $\eta$ function} 
\begin{equation}\label{dedekind}
\eta(\tau)=q^{\frac{1}{12}}\prod_{n=1}^{\infty}(1-q^{2n}),\quad q=e^{i\pi \tau}  
\end{equation}
is a function  satisfying the following relations
\begin{align}\label{rulesdede}
\eta(\tau+1)=&e^{\frac{i\pi}{12}}\eta(\tau),\quad \eta(-1/\tau)=\sqrt{\tau/i}\,\eta(\tau).  
\end{align}

\medskip
\noindent $\bullet$ The {\it theta function} is the exponentially convergent series
\begin{align}\label{theta}
\vartheta_1(z,\tau)=&-i\sum_{n\in\Z}(-1)^n q^{(n+\frac{1}{2})^2}e^{(2n+1)\pi i z},\quad q=e^{i\pi \tau}
\end{align}
defined for $z\in\C$ and $\tau\in\H$. It admits the product representation  
\begin{equation}\label{releteta}
\vartheta_1(z,\tau)=-iq^{\frac{1}{6}}e^{\pi i z}\eta(\tau)\prod_{m=1}^\infty(1-q^{2m}e^{2\pi i z})(1-q^{2m-2}e^{-2\pi i z}).
\end{equation}


\section{Proof of Proposition \ref{shapegreen}}\label{proofshape}
 
Let us set for $s>0$
\begin{equation}
F(s)= \sum_{(n,m)\not=(0,0)}\frac{1}{|n\tau-m|^{\frac{s}{2}}} e^{2\pi i nx_1+2\pi im x_2}.
\end{equation}
Recalling the following relation for $a>0$
$$\int_0^\infty t^{\frac{s}{2}-1}e^{-ta}\,dt=a^{-\frac{s}{2}}\Gamma(\frac{s}{2}),$$
we deduce
$$\Gamma(\frac{s}{2})F(s)= \sum_{(n,m)\not=(0,0)}\int_0^\infty t^{\frac{s}{2}-1}e^{2\pi i nx_1+2\pi im x_2-t|n\tau-m|^2}\,dt.$$
Recall now the Poisson summation formula
\begin{lemma}\label{poisson}
Let $f:\R\to \C$ such that $|f(x)|\leq \frac{C}{(1+|x|)^\alpha}$ for some $C>0$ and $\alpha>1$. Then 
$$\sum_{m}f(x+m)=\sum_{m}\hat{f}(m)e^{2\pi i mx}$$
where $\hat{f}(u)=\int_{\R} e^{-2\pi i u x}f(x)\,dx$.
\end{lemma}

By setting
$$f_n(x)=\sqrt{\frac{ \pi}{t}}e^{- \frac{(2\pi x-i2t n{\rm Re}(\tau))^2}{4 t}-t n^2|\tau|^2},$$
one can check that $\hat{f}_n(m)=e^{-t|n \tau-m|^2}$. By applying Lemma \ref{poisson} for each fixed $n$ and then summing over $n$, we deduce
\begin{align*}
\Gamma(\frac{s}{2})F(s)=&\int_0^\infty t^{\frac{s}{2}-1}\Big(\sum_{m\not=0}e^{2\pi im x_2-t|m|^2}+\sum_{n\not=0,m}f_n(x_2+m)e^{2 i\pi nx_1}\Big)\,dt\\
=&\int_0^\infty t^{\frac{s}{2}-1}\Big(\sum_{m\not=0}e^{2\pi im x_2-t|m|^2}\Big)\,dt+\sqrt{  \pi }\sum_{n\not=0,m} \Big(\int_0^\infty t^{\frac{s-3}{2}}e^{- \frac{(2\pi x_2+2\pi m-i2t n{\rm Re}(\tau))^2}{4 t}-t n^2|\tau|^2}\,dt\Big)e^{2 i\pi nx_1}\\
=&\int_0^\infty t^{\frac{s}{2}-1}\Big(\sum_{m\not=0}e^{2\pi im x_2-t|m|^2}\Big)\,dt\\
&+\sqrt{  \pi }\sum_{n\not=0,m} \Big(\int_0^\infty t^{\frac{s-3}{2}}e^{-\frac{ \pi^2(x_2+m)^2}{t} -t n^2{\rm Im}(\tau)^2}\,dt\Big)e^{2 i\pi nx_1+2i\pi n{\rm Re}(\tau)(x_2+m)}.
\end{align*}
Now recall the Hobson's representation for modified Bessel function of the second kind ($\beta,k,q>0$)
\begin{equation}
\int_0^\infty t^{\beta-1}e^{-k^2 t-\frac{q^2\pi^2}{t}}\,dt=2(q\pi/k)^\beta K_\beta(2\pi k q).
\end{equation}
We use this formula to get
\begin{align*}
\Gamma(\frac{s}{2})F(s)= &  \sum_{m\not=0}\frac{e^{2\pi im x_2}}{|m|^s} +2\sqrt{  \pi }\sum_{n\not=0,m} \Big(\frac{\pi |x_2+m|}{|n|{\rm Im}(\tau)}\Big)^{\frac{s-1}{2}}K_{\frac{s-1}{2}}(2\pi |n| {\rm Im}(\tau)|x_2+m|) e^{2 i\pi nx_1+2i\pi n{\rm Re}(\tau)(x_2+m)}
\end{align*}
Now we take $s=2$ and use the exact expression $K_{\frac{1}{2}}(x)=\sqrt{\frac{\pi}{2x}}e^{-x}$ for $x>0$ to obtain 
\begin{align*}
 F(2)= &  \sum_{m\not=0}\frac{e^{2\pi im x_2}}{|m|^2} + \pi\sum_{n\not=0,m} \frac{1 }{|n|{\rm Im}(\tau)}e^{-2\pi |n| {\rm Im}(\tau)|x_2+m| }e^{2 i\pi nx_1+2i\pi n{\rm Re}(\tau)(x_2+m)}.
\end{align*}
Now we take the sum $n\not=0,m=0$ out of the double sum and then combine the positive/negative values of $m,n$ in the same sum to get
\begin{align*}
 F(2)= &  2\sum_{m\geq 1}\frac{\cos(2\pi m x_2)}{m^2}\\
 & + \frac{\pi }{ {\rm Im}(\tau)}\sum_{n\geq 1} \frac{1 }{n}\Big(e^{-2\pi n {\rm Im}(\tau) x_2 +2 i\pi nx_1+2i\pi n{\rm Re}(\tau)x_2}+e^{-2\pi n {\rm Im}(\tau) x_2 -2 i\pi nx_1-2i\pi n{\rm Re}(\tau)x_2}\Big)\\
 &+\frac{\pi }{ {\rm Im}(\tau)}\sum_{n,m\geq 1}\frac{1}{n}\big(q^{2nm}(e^{2\pi i z})^n+\bar{q}^{2nm}(\overline{e^{-2\pi i z}})^{n}+\bar{q}^{2nm} (\overline{e^{2\pi i z}})^{n}+q^{2nm}(e^{-2\pi i z})^n\big)
\end{align*}
where we have set $q=e^{i\pi \tau}$ and $z=x_1+\tau x_2$. The first sum is standard Fourier series 
$$ \sum_{m\geq 1}\frac{\cos(2\pi m x_2)}{m^2}=\pi^2(x^2-x+\frac{1}{6}).$$ For the remaining sums, we perform the sums over $n$ by using the relation $-\ln(1-x)=\sum_{n\geq 1}\frac{x^n}{n}$ for $x\in\C$ with $|x|<1$. We obtain
\begin{align*}
 F(2)= &  2\pi^2(x_2^2-x_2+\frac{1}{6}) - \frac{2\pi }{ {\rm Im}(\tau)}\ln |1-e^{2i\pi z }| -\frac{2\pi }{ {\rm Im}(\tau)}\ln \prod_{m\geq 1}|1-q^{2m}e^{2\pi i z}||1-q^{2m}e^{-2\pi i z}| .
 \end{align*}
By using the relation \eqref{releteta}, we deduce
$$F(2)=2\pi^2x_2^2-\frac{2\pi}{ {\rm Im}(\tau)}\ln\Big|\frac{\vartheta_1(z,\tau)}{\eta(\tau)}\Big|\quad\text{with } z=p_\tau(x),$$
  which completes the proof.\qed

%
 

\hspace{10 cm}

 \end{document}